\documentclass[a4paper,12pt,final]{amsart}
\usepackage[utf8]{inputenc}
\usepackage[T1]{fontenc}
\usepackage[UKenglish]{babel}
\usepackage[a4paper,margin=28mm]{geometry}
\usepackage{verbatim}

\allowdisplaybreaks[4]
\usepackage{times}
\usepackage{dsfont,mathrsfs}
\usepackage{amsmath}
\usepackage{amsthm}
\usepackage{amssymb}
\usepackage{amsfonts}
\usepackage{latexsym}
\usepackage{booktabs}

\usepackage[colorlinks,linkcolor=red,pagebackref]{hyperref}
\usepackage{xcolor}

\usepackage{ulem}

\newtheorem{theorem}{Theorem}[section]
\newtheorem{lemma}[theorem]{Lemma}

\newtheorem{assumption}[theorem]{Assumption}

\theoremstyle{definition}

\newtheorem{remark}[theorem]{Remark}
\numberwithin{equation}{section}
\renewcommand{\labelenumi}{\roman{enumi})}
\renewcommand\theenumi\labelenumi
\renewcommand{\leq}{\leqslant}
\renewcommand{\le}{\leqslant}
\renewcommand{\geq}{\geqslant}
\renewcommand{\ge}{\geqslant}

\newcommand{\tl}{\tilde}
\newcommand{\Be}{\begin{equation}}
\newcommand{\Ees}{\end{equation*}}
\newcommand{\Bes}{\begin{equation*}}
\newcommand{\Ee}{\end{equation}}

\newcommand{\R}{\mathbb{R}}
\newcommand{\E}{\mathbb{E}}

\newcommand{\N}{\mathbb{N}}
\newcommand{\mcl}{\mathcal}

\newcommand{\dif}{\mathrm{d}}

\usepackage{marginnote}%this avoids floating marginal notes
\usepackage[notref,notcite]{showkeys}

\begin{document}

\title[Approximation of the invariant measure for stable SDE by EM]
{Approximation of the invariant measure for stable SDE by the Euler-Maruyama scheme with decreasing step-sizes}

\author[P.~Chen]{Peng Chen}
\address[P.~Chen]{School of Mathematics, Nanjing University of Aeronautics and
Astronautics, Nanjing 211106, China}
\email{chenpengmath@nuaa.edu.cn}

\author[X. Jin]{Xinghu Jin}
\address[X. ~Jin]{School of Mathematics, Hefei University of Technology,
Hefei, Anhui, China}
\email{xinghujin@hfut.edu.cn}

\author[Y. Xiao]{Yimin Xiao}
\address[Y. Xiao]{%
Department of Statistics and Probability, Michigan State University, 619 Red Cedar Road
East Lansing, MI 48824, United States}
\email{xiao@stt.msu.edu }

\author[L.~Xu]{Lihu Xu}
\address[L.~Xu]{1. Department of Mathematics, Faculty of Science and Technology, University of Macau,
Macau S.A.R., China. 2. Zhuhai UM Science \& Technology Research Institute, Zhuhai, China}
\email{lihuxu@um.edu.mo}

 \subjclass[2010]{}

 \begin{abstract}
Let  $(X_t)_{t \ge 0}$ be the solution of the stochastic differential equation
\[
\dif X_t = b(X_t) \dif t+A \dif Z_t, \quad X_{0}=x,
\]
where  $b: \R^d \rightarrow \R^d$ is a Lipschitz function, $A \in \R^{d \times d}$ is a positive definite matrix,
$(Z_t)_{t\geq 0}$ is a $d$-dimensional rotationally invariant $\alpha$-stable L\'evy process
with $\alpha \in (1,2)$ and $x\in\mathbb{R}^{d}$.
We use two Euler-Maruyama schemes with decreasing step sizes $\Gamma = (\gamma_n)_{n\in \mathbb{N}}$
to approximate the invariant measure of $(X_t)_{t \ge 0}$: one with i.i.d. $\alpha$-stable distributed random
variables as its innovations and the other with i.i.d. Pareto distributed random variables as its innovations. We study
the convergence rate of these two approximation schemes in the Wasserstein-1 distance. For the first scheme,
when the function $b$ is Lipschitz and satisfies a certain
dissipation condition, we show that the convergence rate is  $\gamma^{1/\alpha}_n$. Under an additional assumption
on the second order directional derivatives of $b$, this convergence rate can be improved to $\gamma^{1+\frac 1
{\alpha}-\frac{1}{\kappa}}_n$ for any $\kappa \in [1,\alpha)$.
For the second scheme, when the function $b$ is twice continuously differentiable, we obtain a convergence rate of
$\gamma^{\frac{2-\alpha}{\alpha}}_n$. We show that  the rate $\gamma^{\frac{2-\alpha}{\alpha}}_n$ is optimal
for the one dimensional stable Ornstein-Uhlenbeck process.

Our theorems indicate that the recent remarkable result about the unadjusted Langevin algorithm with additive
innovations in \cite{Pages2020Unajusted} can be extended to the SDEs driven by an $\alpha$-stable L\'evy process
and the corresponding convergence rate has a similar behaviour. Compared with the result in \cite{Chen2022stableEM},
we have relaxed the second order differentiability condition to the Lipschitz condition for the first scheme.
\end{abstract}

\maketitle

% \tableofcontents\thispagestyle{plain}

\section{Introduction and Main results}
%\subsection{Setting}
We consider the following stochastic differential equation (SDE):
\begin{eqnarray}\label{e:SDE}
\dif X_t &=& b(X_t) \dif t+A \dif Z_t, \quad X_{0}=x,
\end{eqnarray}
where the function $b: \R^d \rightarrow \R^d$ satisfies {\bf Assumption \ref{assump-1}} below, $A \in \R^{d \times d}$ is
a positive definite matrix, $(Z_t)_{t\geq 0}$ is a $d$-dimensional rotationally invariant $\alpha$-stable L\'evy  process
with $\alpha \in (1,2)$ and $x\in\mathbb{R}^{d}$.

Under {\bf Assumption \ref{assump-1}} below, \eqref{e:SDE} admits a unique invariant measure $\nu$ (see, e.g., \cite[Corollary 1.4]{Liang2020Gradient}).
%\footnote{\textcolor{blue} {Can you add a reference for this? I see
%that Proposition 1.5 in \cite{Chen2022stableEM} proves a similar result under a stronger condition on $b$.}}
In order to approximate the invariant measure $\nu$, we consider the following Euler-Maruyama (E-M) scheme:
let $Y_{0}=x$,
\begin{eqnarray}\label{e:EM-1}
Y_{t_{n+1}} &=& Y_{t_n} + \gamma_{n+1} b(Y_{t_n}) +A \gamma^{1/\alpha}_n \zeta_{n+1},
\end{eqnarray}
where $t_{0}=0$, $t_{n}=\sum_{i=1}^n \gamma_{i}$, $\{\gamma_{n}\}_{n\in \mathbb{N}}$ is a decreasing sequence of positive
numbers satisfying $\sum_{i=1}^{\infty}\gamma_{i}=\infty$, and $\{\zeta_n\}_{n \in \mathbb N}$ are i.i.d. $d$-dimensional rotationally
invariant $\alpha$-stable random vectors.
%with  standard symmetric $\alpha$-stable distribution.
It is easy to see that the random variables $\{Y_{t_n}\}_{n \ge 1}$ obtained by
\begin{eqnarray}\label{e:EM}
Y_{t_{n+1}} &=& Y_{t_n} + \gamma_{n+1} b(Y_{t_n}) +A (Z_{t_{n+1}}-Z_{t_n}),
\end{eqnarray}
have the same distribution as those obtained by \eqref{e:EM-1}.

Due to the complexity of sampling multi-dimensional stable distributions, see detailed discussion in \cite{Chambers1976} and \cite[Section 1.9]{Nolan2020}, we also use the Pareto distribution for the innovations in the EM scheme \eqref{e:EM2} below. Let $\tl{Z}_1,\tl{Z}_2, \cdots $
be a sequence of i.i.d. $d$-dimensional Pareto distributed random vectors, that is,
\begin{eqnarray}\label{e:Pareto}
\tl{Z}_1 \sim p(z) \ \ = \ \ \frac{\alpha}{\sigma_{d-1} |z|^{\alpha+d}} 1_{(1,\infty)}(|z|).
\end{eqnarray}
Then, we are able to approximate the SDE (\ref{e:SDE}) by the following iterative scheme:
\begin{eqnarray}\label{e:EM2}
\tilde{Y}_{t_{n+1}} &=& \tilde{Y}_{t_n} + \gamma_{n+1} b(\tilde{Y}_{t_n}) + \frac{\gamma_{n+1}^{\frac{1}{\alpha}}}{\beta} A\tl{Z}_{n+1},
\end{eqnarray}
where $\beta=\left(\frac{\alpha}{\sigma_{d-1}d_{\alpha}}\right)^{1/\alpha}$
with $\sigma_{d-1}=2\pi^{\frac{d}{2}}/\Gamma(\frac{d}{2})$ and $d_{\alpha} = \alpha 2^{\alpha-1}\pi^{-d/2}\Gamma(\frac{d+\alpha}{2})/\Gamma(1-\frac{\alpha}{2})$.

\subsection{Motivations and related work}
There have been plenty of works on the EM scheme for SDEs, see, for example, \cite{Bally1996The,Bao2019convergence,Dareiotis2020On,Janicki1996Approximation,
Li2023Strong,Jacod2004The,Li2018The,Marion2022Convergence,Sanz2021Wasserstein, Shao2018Invariant,Tarami2018Convergence,Yuan2008a}
and the references therein for more information. In the past decade, the EM schemes for the SDEs like \eqref{e:SDE},
have been intensively studied by the authors of \cite{BDG22,Has13,Huang2018The,Kuhn2019Strong,LZ23,
Mimica2018Markov,Pamen2017Strong,Panloup2008Recursive,
Protter1997The}. Most of them have only considered the case where the step size of the EM scheme is constant and
proved the convergence rate in a finite time interval $[0,T]$ for a constant $T>0$. As $T \rightarrow \infty$, the rate
usually blows up. To the best of our knowledge, \cite{Chen2022stableEM} is the first paper that has studied the long term
behaviour of the EM scheme of the SDE driven by an $\alpha$-stable process, but the step size therein is constant. However,
in most of real applications, the step sizes (or learning rate) decrease as the iteration goes on. What's more, \cite{Chen2022stableEM}
adds a strong condition on the drift function $b$ in the first EM scheme. One of our motivations for the research in this paper
is to improve the results in \cite{Chen2022stableEM}.

Another motivation for this paper is from the recent paper \cite{Pages2020Unajusted} by Pag\`{e}s and Panloup, who
considered the convergence rate of an unadjusted Langevin algorithm with decreasing step sizes. Our EM schemes have
a similar form as the unadjusted Langevin algorithm with additive innovation in \cite{Pages2020Unajusted}, only replacing
the standard normal random variables therein by the Pareto random vectors or $\alpha$-stable random vectors. A remarkable
result in \cite{Pages2020Unajusted} is that the unadjusted Langevin algorithm converges to the targeted measure at a rate
only depending on the step sizes $\gamma_n$ at the terminal iterate $n$. We shall prove in this paper that the EM scheme
\eqref{e:EM2} has the same phenomenon in the Wasserstein-1 distance setting.

%\footnote{\textcolor{blue}{A phrase that is similar to "There have been plenty of works on
%the EM scheme for SDEs" has appeared at the beginning of this paragraph and the references have some overlap. Would it
%be better to combine them?}}

Due to the significance of high dimensional sampling problems arising in machine learning, there has recently been a
surge in studying the convergence rates of the EM schemes to their corresponding SDEs as time tends to infinity. Fang et al.
\cite{Fang2019Multivariate} used the EM scheme to approximate the invariant measure for an SDE driven by Brownian
motion with help of the Malliavin calculus and Stein's method. Besides the recent developments on the
convergence analysis of unadjusted Langevin algorithms \cite{D17, D-M17, D-M19}, heavy tailed stochastic
algorithms have recently attracted a lot of attention, see for instance \cite{LWS2022NIPS,NDHR2021COLT,NSGR2019NIPS}.
In particular, by using the convergence result in \cite{Chen2022stableEM},  the heavy tailed stochastic gradient descent algorithms
are delicately analyzed in \cite{RBGZS23ICML,Raj2022Algorithmic}. We expect that the more general results in this paper will be applicable for studying heavy tailed stochastic algorithms.

%Pag\`{e}s and Panloup \cite{Pages2020Unajusted} used the EM scheme with decreasing step size to study the approximation of invariant measure of SDE driven by Brownian motion with help of the Malliavin calculus. In this paper, we shall use this kind of method to discuss the approximation of the invariant measure for SDE driven by symmetric $\alpha$-stable processes.
%Recently, . There are also some references concerning on the computation of the invariant measure for SDEs, see \cite{Benaim2017Ergodicity,Fang2020Adaptive,Panloup2008Computation,Panloup2008Recursive} and references therein, most of them are asymptotic type. Strong convergence for EM scheme can be found in \cite{Chen2019Numerical,Hoang2014Numerical,Kuhn2019Strong,Yin2005Numerical} and references therein.

%We concern on the symmetric $\alpha$-stable process in this paper rather than the Brownian motion, the approximation of stable law can be found in \cite{Chen2019Approximation, Jin2020A,Xu2019Approximation} and references therein. \cite{Huang2018The,Kuhn2019Strong,Pamen2017Strong,Panloup2008Recursive,Panloup2008Computation,Protter1997The} and references therein concerned on the approximation for invariant measure for SDEs driven by L\'{e}vy processes.

\subsection{Notations}
We denote by $|.|$ the Euclidean norm of $\mathbb R^d$, i.e., for $x \in \mathbb R^d$, $|x|=\sqrt{x^2_1+...+x^2_d}$.
Given any two probability measures $\nu_1$ and $\nu_2$ on $\R^d$, their $L^p$-Wasserstein distance $W_p$
with $p\geq 1$ is defined as
\begin{eqnarray}\label{e:Wp}
W_p(\nu_1,\nu_2) &=& \inf_{ \Pi \in \mathcal{C}(\nu_1,\nu_2) } \left( \int_{\R^{2d}}  |x-y|^p  \dif \Pi(\dif x,\dif y)  \right)^{\frac{1}{p}},
\end{eqnarray}
where $\mathcal{C}(\nu_1,\nu_2)$ is the set of probability distributions on $\R^d \times \R^d$ with the marginal
measures $\nu_1$ and $\nu_2$ respectively. By using the Kantorovich duality (see \cite[Theorem 5.10]{Villani2009Optimal}),
$W_1$ can also be defined as
\begin{eqnarray}\label{e:W1}
W_1(\nu_1,\nu_2) &=& \sup_{ h\in {\rm Lip}(1)} \left(   \int_{\R^d} h(x) \nu_1(\dif x)  -  \int_{\R^d} h(x) \nu_2(\dif x)   \right),
\end{eqnarray}
where ${\rm Lip}(1)$ is the set of Lipschitz functions with Lipschitz constant $1$, that is, ${\rm Lip}(1)
= \{ f: \frac{|f(x)-f(y)|}{|x-y|} \leq 1, \ \ \forall x, y \in \R^d   \}.$

For matrices $A=(A_{ij})_{d\times d}$ and $B=(B_{ij})_{d\times d}$, define $ \langle A,B\rangle_{\rm HS}=\sum_{i,j=1}^d A_{ij} B_{ij}$
and the Hilbert-Schmidt norm is $\| A\|_{\rm HS}=\sqrt{ \langle A,A\rangle_{\rm HS} }$.

For $k \in \N_0:=\N \cup\{0\}$, denote by $\mcl C^k\left(\mathbb{R}^d, \mathbb{R}\right)$ the set of functions which have
$k$-th order continuous derivatives. In particular, $\mcl C^0\left(\mathbb{R}^d, \mathbb{R}\right)= \mcl C\left(\mathbb{R}^d, \mathbb{R}\right)$
denotes the set of continuous functions. For $f \in \mcl C^2\left(\mathbb{R}^d, \mathbb{R}\right)$ and $v_1, v_2, x \in \mathbb{R}^d$,
the directional derivatives $\nabla_{v_1} f(x)$ and $\nabla_{v_2} \nabla_{v_1} f(x)$ are defined by
$$
\nabla_{v_1} f(x):=\lim _{\epsilon \rightarrow 0} \frac{f\left(x+\epsilon v_1\right)-f(x)}{\epsilon},
$$
$$
\nabla_{v_2} \nabla_{v_1} f(x):=\lim _{\epsilon \rightarrow 0} \frac{\nabla_{v_1} f\left(x+\epsilon v_2\right)-\nabla_{v_1} f(x)}{\epsilon}.
$$
The operator norms of $\nabla^2 f(x) \in \mathbb{R}^{d \times d}$ is defined by
$$
\left\|\nabla^2 f(x)\right\|_{\mathrm{op}}:=\sup \left\{\left|\nabla_{v_2} \nabla_{v_1} f(x)\right| ; v_1, v_2 \in \mathbb{R}^d,
\left|v_1\right|=\left|v_2\right|=1\right\}.
$$
Denote $\| f \|_{\infty} = \sup_{x\in \R^d}|f(x)|$ for a bounded continuous function $f$. For any $f\in \mathcal{C}^2(\R^d,\R)$,
let $\nabla f$ and $\nabla^2 f $ be the gradient and Hessian matrix of $f$, and denote $\| \nabla f \|_{\infty} = \sup_{x\in \R^d} |\nabla f(x)|$
and $\| \nabla^2 f \|_{{\rm HS},\infty}= \sup_{x\in \R^d} \| \nabla^2 f(x) \|_{{\rm HS}}$ if they both exist.

For  $k \in \N_0$, further denote by $\mcl C^k_{lin}(\R^d,\R)$ the set of functions in $\mcl C^k\left(\mathbb{R}^d, \mathbb{R}\right)$
that have a linear growth, i.e.
$$\mcl C^k_{lin}(\R^d,\R)=\left\{ f \in \mcl C^k(\R^d,\R):  \ \sup_{x \in \R^d} \frac{|f(x)|}{1+|x|}<\infty\right\}.$$
As $k=0$, we write it as $\mcl C_{lin}(\R^d,\R)$.

For any random variable $X$, denote its law by $\mathcal{L}(X)$ and for any $r_{1},r_{2}\in\mathbb{R}$, denote
$r_{1}\vee r_{2}=\max\{r_{1},r_{2}\}$.

%\begin{assumption}\label{assump-1}
%(i) Let $b\in \mathcal{C}^2(\R^d,\R^d)$ and there exist $\theta_1,  L>0$ and $\theta_{2}, K\geq 0$ such that
%\begin{eqnarray}\label{e:b-1}
%\langle x-y, b(x)-b(y) \rangle &\leq& -\theta_1 |x-y|^2+K, \quad \forall x,y \in \R^d,
%\end{eqnarray}
%and
%\begin{eqnarray}\label{e:b-2}
%|\nabla_{u_1} b(x)|  \ \ \leq \ \   L |u_1|, \quad
%|\nabla_{u_2} \nabla_{u_1} b(x)|  \ \ \leq \ \ \theta_{2} |u_1| |u_2|,
%\quad   \forall x, v_1, v_2 \in \R^d.
%\end{eqnarray}
%\end{assumption}

\subsection{Assumptions and main result}

We put the following assumptions on the drift $b$.
\begin{assumption}\label{assump-1} Let $b$ be a Lipschitz function, i.e., there exists a constant $L>0$ such that
\begin{eqnarray}\label{e:b-1}
|b(x)-b(y)| &\leq& L |x-y| \quad \quad \forall \ x,y \in \R^d.
\end{eqnarray}
Moreover, there exist constants $\theta_1$ and $K\geq 0$ such that
\begin{eqnarray}\label{e:b-1b}
\langle x-y, b(x)-b(y) \rangle &\leq& -\theta_1 |x-y|^2+K, \quad \forall x,y \in \R^d.
\end{eqnarray}
\end{assumption}
%and
%\begin{eqnarray}\label{e:b-2}
%|\nabla_{u_1} b(x)|  \ \ \leq \ \   L |u_1|, \quad
%|\nabla_{u_2} \nabla_{u_1} b(x)|  \ \ \leq \ \ \theta_{2} |u_1| |u_2|,
%\quad   \forall x, u_1, u_2 \in \R^d.
%\end{eqnarray}
%where $\nabla_{u_1}$ and  $\nabla_{u_1} \nabla_{u_2}$ are directional derivatives defined below.

For the sequence $(\gamma_{n})_{n\in \mathbb{N}}$ of step sizes, we consider the following assumption:
\begin{assumption}\label{stepsize}
(i) $(\gamma_n)_{n\in \mathbb{N}}$ is a decreasing sequence, $\lim\limits_{n\to \infty} \gamma_n=0$ and
$\sum_{n=1}^{\infty}\gamma_{n}=\infty$.

(ii) $(\gamma_n)_{n\in \mathbb{N}}$ satisfies the following inequality:
\begin{eqnarray}\label{omega}
\omega := \limsup_{k\to \infty} \frac{\gamma_k^{\theta}-\gamma_{k+1}^{\theta}}{\gamma_{k+1}^{1+\theta}} <+\infty
\end{eqnarray}
for some $\theta\in(0,1]$.
\end{assumption}

Now we are at the position to state our main theorems. We first give an error bound between the processes
$(X_t)_{t\geq 0}$ and $(Y_{t_n})_{n\in \mathbb{N}_0}$. In the following, $\mathcal{L}(X_{t}^x)$ denotes the
law of $X_t$ starting from $x$ at time $t = 0$. The notation $\mathcal{L}(Y_{t_n}^x)$ is similarly defined.

\begin{theorem}\label{thm:LXLY12}
Let $(X_t)_{t\geq 0}$ and $(Y_{t_n})_{n\in \mathbb{N}_0}$ be governed by \eqref{e:SDE} and \eqref{e:EM} respectively.
Under Assumption \ref{assump-1}, suppose that Assumption \ref{stepsize} holds with $\theta=\frac{1}{\alpha}$ and
$\omega<\rho=d^{2\alpha-4}2^{-d} e^{-2^{d} d^{4-2\alpha}}$. Then, there exists a positive constant $C$ independent of
$d$ and $n$ such that for all $x \in \R^d$ and $t_n>1$,
\begin{align}\label{Liptheorem}
W_1(\mathcal{L}(X_{t_n}^x), \mathcal{L}(Y_{t_n}^x))
\leq\frac{C}{\rho-\omega}
e^{\frac{\rho-\omega}{2}\gamma_{1}}
e^{d^{4-2\alpha} 2^{d}} d^{5}
\, \gamma_{n}^{\frac{1}{\alpha}} (1+|x|),
\end{align}
and
\begin{align}\label{LiptheoremB}
W_1(\nu,\mathcal{L}(Y_{t_n}^x))
\leq \frac{C}{\rho-\omega}
e^{\frac{\rho-\omega}{2}\gamma_1} e^{d^{4-2\alpha} 2^{d}} d^{5}\,
\gamma_{n}^{\frac{1}{\alpha}} (1+|x|).
\end{align}
Here and below, $\nu$ denotes the invariant measure of $(X_t)_{t\geq 0}$.
\end{theorem}

\begin{remark}  As an example of the decreasing step sizes $(\gamma_n)_{n\in \mathbb{N}}$, we take
\begin{eqnarray*}
\gamma_k &=& \frac{c}{\rho k}, \ \ \ \forall k\in \mathbb{N},
\end{eqnarray*}
where $\rho =d^{2\alpha-4}2^{-d}
e^{-2^{d} d^{4-2\alpha}}$ as in Theorem \ref{thm:LXLY12}  and $c>1$ is a constant.
Let $\theta=\frac{1}{\alpha}$ in inequality  \eqref{omega}, a straightforward calculation yields
\begin{eqnarray*}
\omega &=& \limsup_{k\to \infty} \frac{\gamma_{k}^{\frac{1}{\alpha}}-\gamma_{k+1}^{\frac{1}{\alpha}}}
{\gamma_{k+1}^{1+\frac{1}{\alpha}}}
\ \ = \ \ \frac{\rho}{\alpha c}
\ \ < \ \ \rho.
\end{eqnarray*}
In this case, the constant $\frac{C}{\rho-\omega}
e^{\frac{\rho-\omega}{2}\gamma_{1}}
e^{d^{4-2\alpha} 2^{d}}d^{5}$ in Theorem \ref{thm:LXLY12} reduces to $Ce^{d^{4-2\alpha}2^{d+1}} d^{5}$.
\end{remark}

We also consider an additional smoothness condition for the coefficient $b$, that is,

\begin{assumption}\label{smooth}
Let $b\in\mathcal{C}^{2}(\mathbb{R}^{d})$ and there exists a constant $\theta_{2}>0$ independent of dimension
$d$ such that
\begin{align*}
\left|\nabla_{u_{2}}\nabla_{u_{1}}b(x)\right|\leq\theta_{2}|u_{1}||u_{2}|, \quad \forall x,u_{1},u_{2}\in\mathbb{R}^{d}.
\end{align*}
\end{assumption}

\begin{remark}
It is easy to verify from Assumption \ref{assump-1} that
\begin{align}\label{Linear}
|b(x)|\leq |b(0)|+L|x|, \ \ \forall x\in \R^d.
\end{align}
Inequality \eqref{e:b-1} and Assumption \ref{smooth} further imply
\begin{align}\label{first-order}
\left|\nabla_{u}b(x)\right|\leq L|u|, \quad \forall x,u\in\mathbb{R}^{d}.
\end{align}
\end{remark}

\begin{theorem}\label{thm:LXLY1}
Let $(X_t)_{t\geq 0}$ and $(Y_{t_n})_{n\in \mathbb{N}_0}$ be governed by \eqref{e:SDE} and \eqref{e:EM} respectively.
Let Assumptions \ref{assump-1} and \ref{smooth} hold and suppose that Assumption \ref{stepsize} holds with $\theta=
1+\frac{1}{\alpha}-\frac{1}{\kappa}$ and $\omega<\rho$, where $\kappa \in[1,\alpha)$ and $\rho=d^{2\alpha-4}2^{-d}
e^{-2^{d} d^{4-2\alpha}}$. Then, there exists a constant $C>0$ independent of $d$ and $n$ such that for all $x \in \R^d$
and $t_n>1$,
\begin{eqnarray}\label{e:w-distance1}
&&W_1(\mathcal{L}(X_{t_n}^x), \mathcal{L}(Y_{t_n}^x))\nonumber\\
&\leq& \frac{C}{\rho-\omega}e^{\frac{\rho-\omega}{2}\gamma_1} e^{d^{4-2\alpha} 2^{d}} d^{5\vee(3\kappa)\vee(4
+\frac{3}{2\kappa}) \vee (\frac{3}{2}+\frac{3}{2\kappa}+\frac{3\kappa}{2} )}\,
\gamma_{n}^{1+\frac{1}{\alpha}-\frac{1}{\kappa}}(1+|x|^{\kappa})
\end{eqnarray}
and
\begin{align}\label{e:w-distance1B}
W_1(\nu,\mathcal{L}(Y_{t_n}^x))\leq \frac{C}{\rho-\omega}e^{\frac{\rho-\omega}{2}\gamma_1} e^{d^{4-2\alpha} 2^{d}}
d^{5\vee(3\kappa)\vee(4+\frac{3}{2\kappa}) \vee (\frac{3}{2}+\frac{3}{2\kappa}+\frac{3\kappa}{2})}\,
\gamma_{n}^{1+\frac{1}{\alpha}-\frac{1}{\kappa}} (1+|x|^{\kappa}).
\end{align}
%where $\nu$ is the invariant measure of $(X_t)_{t\geq 0}$.
\end{theorem}

For the EM scheme with Pareto distributed noise, we have the following convergence theorem.
\begin{theorem}\label{thm:LXLY}
Let $(X_t)_{t\geq 0}$ and $(\tilde{Y}_{t_n})_{n\in \mathbb{N}_0}$ be governed by \eqref{e:SDE} and \eqref{e:EM2}
respectively. Let Assumptions \ref{assump-1} and \ref{smooth} hold and suppose that Assumption \ref{stepsize} holds
with $\theta=\frac{1}{\alpha}$ and $\omega<\rho= d^{2\alpha-4}2^{-d} e^{-2^{d} d^{4-2\alpha}}$. Then, there exists
a positive constant $C$ independent of $d$ and $n$ such that for  all $x \in \R^d$ and $t_n>1$,
\begin{equation}\label{e:w-distance}
W_1\big(\mathcal{L}(X_{t_n}^x), \mathcal{L} (\tilde{Y}_{t_n}^x)\big)
\leq \frac{C} {\rho-\omega} e^{\frac{\rho-\omega} {2}\gamma_1} e^{d^{4 - 2\alpha} 2^{d}} d^{5}\,
\gamma_{n}^{\frac{2-\alpha} {\alpha} }(1+|x|)
\end{equation}
and
\begin{align}\label{e:w-distanceB}
W_1(\nu,\mathcal{L}(\tilde{Y}_{t_n}^x))\leq \frac{C}{\rho-\omega}e^{\frac{\rho-\omega}{2}\gamma_1} e^{d^{4-2\alpha} 2^{d}}d^{5}
 \gamma_{n}^{\frac{2-\alpha} {\alpha}} (1+|x|).
\end{align}
\end{theorem}

\begin{remark}
\begin{itemize}
\item[(i)] For any $0<\theta\leq \delta\leq1$, by the fact
\begin{align*}
\frac{\gamma_{n}^{\theta}-\gamma_{n+1}^{\theta}}{\gamma_{n+1}^{1+\theta}}\leq\frac{\gamma_{n}^{\delta}-\gamma_{n+1}^{\delta}}
{\gamma_{n+1}^{1+\delta}},
\end{align*}
we remark that the condition on the step sequence in Theorem \ref{thm:LXLY1} is slightly stronger than that in Theorems \ref{thm:LXLY12}
and \ref{thm:LXLY}.
\item[(ii)] In Theorem \ref{thm:LXLY}, in order to ensure the convergence of the series in the proof, we need Assumption \ref{stepsize}
with $\theta=\frac{1}{\alpha}$, which is slightly stronger than the condition $\theta=\frac{2-\alpha}{\alpha}$.
\item[(iii)]The rate $\gamma_{n}^{\frac{2}{\alpha}-1}$ in Theorem \ref{thm:LXLY} is optimal for the one-dimensional stable
Ornstein-Uhlenbeck (OU) process, see Lemma \ref{lower bound} in Appendix \ref{seclow}.
\end{itemize}
\end{remark}

The rest of this paper is organized as follows. %We end this section by giving the organization of this paper.
We will give the proofs of our main results, Theorems \ref{thm:LXLY12}, \ref{thm:LXLY1}, and \ref{thm:LXLY}, in Section \ref{sec:proof-thm}.
In Section \ref{APP:ss}, we will prove the crucial lemmas stated in Section \ref{sec:proof-thm}. Finally, in Appendix \ref{seclow}, we show that
the optimal convergence rate $\gamma_{n}^{\frac{2}{\alpha}-1}$ of the EM scheme with the Pareto distribution is reached for the
stable Ornstein-Uhlenbeck process on $\mathbb{R}$.

\section{Proof of Theorems \ref{thm:LXLY12}, \ref{thm:LXLY1} and \ref{thm:LXLY}} \label{sec:proof-thm}
For any $t_j>t_i \ge 0$ and $f \in \mcl C_{lin}(\R^d,\R)$, define
\begin{eqnarray*}
P_{t_i,t_j} f(x)&=&\E[f(X_{t_j})|X_{t_i}=x], \ \ \ \forall \ x \in \R^d, \\
Q_{t_i,t_j} f(x)&=&\E[f(Y_{t_j})|Y_{t_i}=x], \ \ \ \forall \ x \in \R^d,
\end{eqnarray*}
and
\begin{eqnarray*}
\tilde{Q}_{t_i,t_j} f(x) &=& \E[f(\tilde{Y}_{t_j})|\tilde{Y}_{t_i}=x], \ \ \ \forall \ x \in \R^d.
\end{eqnarray*}
Then, using the discrete version of Duhamel principle, we split $P_{0,t_n} h(x) - Q_{0,t_n} h(x)$ as below:
\begin{eqnarray}\label{Duhamel}
P_{0,t_n} h(x) - Q_{0,t_n} h(x)
&=& \sum_{i=1}^{n} Q_{0,t_{i-1}} \circ ( P_{t_{i-1},t_{i}}  - Q_{t_{i-1},t_{i}}  ) \circ P_{t_{i},t_n} h(x) \\
&=& \sum_{i=1}^{n^*} Q_{0,t_{i-1}} \circ ( P_{t_{i-1},t_{i}}  - Q_{t_{i-1},t_{i}}  ) \circ P_{t_{i},t_n} h(x) \nonumber  \\
&&+ \sum_{i=n^*+1}^{n-1} Q_{0,t_{i-1}} \circ ( P_{t_{i-1},t_{i}}  - Q_{t_{i-1},t_{i}}  ) \circ P_{t_{i},t_n} h(x) \nonumber \\
&&+ Q_{0,t_{n-1}} \circ ( P_{t_{n-1},t_{n}}  - Q_{t_{n-1},t_{n}}  )  h(x), \nonumber
\end{eqnarray}
where $n^*=\max\{i:t_n-t_i>1\}$.

%In order to bound the terms in (\ref{Duhamel})
For proving Theorems \ref{thm:LXLY12}, \ref{thm:LXLY1} and \ref{thm:LXLY}, the exponential contraction
property and the regularity of $P_{t_i, t_n}$ play a crucial role. We state them as lemmas in the next section.
%which be addressed below,
% Propositions  \ref{pro:Pt} and \ref{pro:P-Qf} are proven in Sections \ref{sec:Malliavin}  and \ref{sec:P-Qf} respectively.
%Moments estimates in Lemmas \ref{lem:X-moment}, \ref{lem:Y-moment} and \ref{lem:X-X-moment} are postponed in Section \ref{sec:ergodic}.

\subsection{Auxiliary propositions and lemmas}

We first give the moment estimates for $(X_t)_{t\geq 0}$, $(Y_{t_n})_{n\in \mathbb{N}_0}$ and $(\tl{Y}_{t_n})_{n\in \mathbb{N}_0}$,
which will be proved in Section \ref{APP:ss}. For simplicity, for any $i\leq j$, we denote by $X_{t_{i},t_{j}}^{x}$ %to emphasize
the value of the process at time $t_j$ %$X_{t_{i},t_{j}}$
while starting from $x$ at the time $t_{i}$. If $t_i=0$, then we denote it as $X_{t_j}^x$. The notations $Y_{t_{i},t_{j}}^{x}$
and $\tilde{Y}_{t_{i},t_{j}}^x$ are similarly defined.
%In addition, for any $\kappa \in[1,\alpha)$, denote
%\begin{eqnarray}\label{e:V}
%V_{\kappa}(x) &=&(1+|x|^2)^{\frac{\kappa}{2}}, \quad x\in \R^d.
%\end{eqnarray}

%Then, the following proposition is from \cite[Proposition 1.5]{Chen2022stableEM}.

%\begin{proposition}\label{lem:X-moment}
%Let Assumption \ref{assump-1} hold and denote by $(X_{t})_{t\geq0}$ the solution to the SDE \eqref{e:SDE}. Then,
%$(X_{t})_{t\geq0}$ admits a unique invariant probability measure $\nu$ such that for $1\leq\beta<\alpha$,
%\begin{align*}
%\sup_{|f|\leq V_{\beta}}\left|\mathbb{E}[f(X_{t}^{x})]-\nu(f)\right|\leq c_{1}V_{\beta}(x)e^{-c_{2}t}, \quad t\geq0,
%\end{align*}
%for some constants $c_{1},c_{2}>0$. In particular, there exists a constant $C>0$ such that
%\begin{eqnarray*} % \label{e:X-moment}
%\E |X_t^x|\leq C(1+|x|^{\beta}), \ \ t\geq 0.
%\end{eqnarray*}
%\end{proposition}

\begin{lemma}\label{lem:Y-moment}
Let $(X_{t})_{t\geq0}$, $(Y_{t_n})_{n\in \mathbb{N}_0}$ and $(\tilde{Y}_{t_n})_{n\in \mathbb{N}_0}$ be the processes
defined by \eqref{e:SDE}, \eqref{e:EM} and (\ref{e:EM2}), respectively. Suppose that Assumption \ref{assump-1} holds,
then for any $1\leq\kappa <\alpha$, there exists a constant  $C>0$ independent of $d$ such that for all $x \in \R^d$,
\begin{eqnarray}
\E |X_{t}^x|^{\kappa} &\leq& Cd^{\frac{5\vee(3\kappa)}{2}}(1+|x|^{\kappa}), \ \  \forall t\geq0, \label{X-moment} \\
\E |Y_{t_n}^x|^{\kappa} &\leq& Cd^{\frac{5\vee(3\kappa)}{2}}(1+|x|^{\kappa}), \ \  \forall n\in \mathbb{N}_0, \label{e:Y-moment1} \\
\E |\tilde{Y}_{t_n}^x|^{\kappa} &\leq& Cd^{\frac{5\vee(3\kappa)}{2}} (1+|x|^{\kappa}) , \ \  \forall n\in \mathbb{N}_0. \label{e:Y-moment}
\end{eqnarray}
\end{lemma}

\begin{lemma}\label{lem:mom}
Let $(X_{t})_{t\geq 0}$ be the solution to the SDE \eqref{e:SDE} and let $(Y_{t_{n}})_{n\in\mathbb{N}_{0}}$ be the
Markov chains defined by \eqref{e:EM}. Suppose that Assumption \ref{assump-1} holds, then for any $1\leq\kappa <\alpha$,
there exists a positive constant $C$ independent of $d$ such that for all $n\in\mathbb{N}$ and $x\in\mathbb{R}^{d}$:
\begin{align}\label{1}
\mathbb{E} |Y_{t_{n-1},t_{n}}^{x}-x|^{\kappa} &\leq Cd^{\kappa} (1+|x|^{\kappa})\gamma_{n}^{\kappa/\alpha},\\
\label{2}
\mathbb{E} |X_{s,t}^{x}-x|^{\kappa}&\leq C d^{\frac{5\vee(3\kappa)}{2}} (1+|x|^{\kappa})(t-s)^{\kappa/\alpha}, \quad \forall 0\leq s\leq t\leq \gamma_{1},\\
\label{3}
\mathbb{E} |X_{t_{n-1},t_{n}}^{x}-Y_{t_{n-1},t_{n}}^{x}|^{\kappa}&\leq Cd^{\frac{5\vee(3\kappa)}{2}}(1+|x|^{\kappa})\gamma_n^{\kappa+\frac{\kappa}{\alpha}}.
\end{align}
\end{lemma}

Based on the above lemmas, according to \cite[Lemma 18]{CNXY23}, we can derive the following bound on the difference
between $P_{t_{n-1},t_{n}}f$ and $Q_{t_{n-1},t_{n}}f$, the proof will be postponed to Section \ref{APP:ss}.
\begin{lemma}\label{pro:P-Qf1}
Under Assumptions \ref{assump-1} and \ref{smooth}, let $f\in \mathcal{C}^2(\R^d,\R)$ satisfying $\| \nabla f \|_{\infty} < \infty$
and $\| \nabla^2 f \|_{{\rm HS}, \infty} < \infty$. Then for $\kappa \in[1,\alpha)$, there exists a constant $C>0$ independent of
$d$ such that for all $k\in \mathbb{N}$ and $x\in\mathbb{R}^{d}$,
\begin{align*}
| (P_{t_{k-1},t_{k}}  - Q_{t_{k-1},t_{k}}  )  f(x) |
\leq& C(1+|x|^{\kappa}) (\| \nabla f \|_{\infty} +  \| \nabla^2 f \|_{{\rm HS},\infty})d^{\frac{5}{2}\vee\frac{3\kappa}{2}\vee(\frac{3}{2}
+\frac{3}{2\kappa})} \gamma_k^{2+\frac{1}{\alpha}-\frac{1}{\kappa}}.
\end{align*}
\end{lemma}

In the second EM scheme, the difference between $P_{t_{n-1},t_{n}}f$ and $\tilde{Q}_{t_{n-1},t_{n}}f$ is estimated in the following lemma.

\begin{lemma}\label{pro:P-Qf}
 Under Assumption \ref{assump-1}, let $f\in \mathcal{C}^2(\R^d,\R)$ satisfying $\| \nabla f \|_{\infty} < \infty$ and $\| \nabla^2 f \|_{{\rm HS}, \infty} $
$ < \infty$. There exists a constant $C>0$ independent of $d$ such that for all $k\in \mathbb{N}$ and $x\in\mathbb{R}^{d}$,
\begin{eqnarray*}
| (P_{t_{k-1},t_{k}}  - \tilde{Q}_{t_{k-1},t_{k}}  )  f(x) |
&\leq& C(1+|x|) (\| \nabla f \|_{\infty} +  \| \nabla^2 f \|_{{\rm HS},\infty}) d^{\frac{5}{2\kappa} \vee (3-\frac{\alpha}{2}) }  \gamma_k^{\frac{2}{\alpha}}.
\end{eqnarray*}
\end{lemma}

We also need the following regularities for $P_{t_i,t_j}h(x)$ with $h\in {\rm Lip}(1)$ and $j\geq i$, which can be proved in a similar way as
in \cite[Lemma 2.1]{Chen2022stableEM}.
\begin{lemma}\label{pro:Pt}
Let $h\in {\rm Lip}(1)$ and let Assumptions \ref{assump-1} and \ref{smooth} hold. Then there exists a constant $C>0$ independent
of $d$ such that for all $0<t_j-t_i \le 1$ and $u_1,u_2 \in \R^d$,
\begin{eqnarray*}
\| \nabla_{u_1} P_{t_i,t_j} h \|_{\infty} &\leq& C|u_1|, \\
\| \nabla_{u_2} \nabla_{u_1} P_{t_i,t_j} h \|_{{\rm HS},\infty} &\leq& C(t_j-t_i)^{-\frac{1}{\alpha}} |u_1| |u_2|.
\end{eqnarray*}
\end{lemma}

The exponential ergodicity of the solution  $(X_t)_{t\geq 0}$ of SDE \eqref{e:SDE} under Assumption \ref{assump-1}
can be obtained by using the same argument as the proof of \cite[Theorem 1.3]{Liang2020Gradient} and the proof will
be postponed to Section \ref{sec:para}.
\begin{lemma}\label{pro:erg-X}
Under Assumption \ref{assump-1}, there exists a constant $C>0$ independent of $d$ such that for all $x,y\in\mathbb{R}^{d}$
and $t>0$,
\begin{eqnarray*}
W_1 (\mathcal{L}(X_{t}^x), \mathcal{L}(X_{t}^y)) &\leq& Ce^{ d^{4-2\alpha} 2^{d} } e^{-d^{2\alpha-4}2^{-d}e^{-2^{d}d^{4-2\alpha}} t}|x-y|.
\end{eqnarray*}
\end{lemma}

In addition, the following lemma about the decreasing steps also plays an important role in proving the main theorems.
This lemma will be proved in Section \ref{APP:ss}.

\begin{lemma}\label{lem:cal-ss}
For some $\theta\in(0,1]$, let $(\gamma_{n})_{n\in \mathbb{N}}$ be a decreasing positive sequence such that
\begin{eqnarray*}
\omega := \limsup_{k\to \infty} \frac{\gamma_k^{\theta}-\gamma_{k+1}^{\theta}}{\gamma_{k+1}^{1+\theta}} <+\infty.
\end{eqnarray*}
Let $\omega<\rho= d^{2\alpha-4}2^{-d}
e^{-2^{d} d^{4-2\alpha}}$ and let $(v_{n})_{n\in \mathbb{N}}$ be the sequence defined by $v_{0}=0$ and, for every $n\geq1$, by
\begin{align*}
v_{n}=\sum_{k=1}^{n}\gamma_{k}^{1+\theta}e^{-\rho(t_{n}-t_{k})}.
\end{align*}
Then,
\begin{eqnarray}
\limsup_{n\rightarrow\infty}\frac{v_{n}}{\gamma_{n}^{\theta}}
&\leq&
\frac{2}{\rho-\omega}e^{\frac{\rho-\omega}{2}\gamma_1}, \label{convergence1} \\
\lim_{n\rightarrow\infty}\frac{e^{-\rho t_{n}}}{\gamma_{n}^{\theta}}
&=& 0. \label{convergence2}
\end{eqnarray}
Assume that $n$ is large enough such that $t_{n}>1$ and denote $n^*=\max\{ i: t_n - t_i >1 \}$, we further have
\begin{align}\label{convergence3}
\sum_{i=n^*+1}^{n-1}(t_n-t_i)^{-\frac{1}{\alpha}} \gamma_i^{1+\theta}
\leq C e^{\rho(1+\gamma_{1})} \gamma_{n}^{\theta}
\end{align}
for some constant $C>0$ independent of $d$.
\end{lemma}

\subsection{Proof of Theorem \ref{thm:LXLY12}}
(i) It follows from \eqref{e:W1} and \eqref{Duhamel} that
\begin{align}\label{framework}
W_1(\mathcal{L}(X_{t_n}^x), \mathcal{L}(Y_{t_n}^x))
=& \sup_{h\in {\rm Lip}(1)}  \{ P_{0,t_n} h(x) - Q_{0,t_n} h(x)   \}\nonumber\\
\leq&  \sum_{i=1}^{n}\sup_{h\in{\rm Lip}(1)}\left|Q_{0,t_{i-1}} \circ ( P_{t_{i-1},t_{i}}  - Q_{t_{i-1},t_{i}}  ) \circ P_{t_{i},t_n} h(x)\right|\nonumber\\
=& \sum_{i=1}^{n-1}\sup_{h\in{\rm Lip}(1)}\left| Q_{0,t_{i-1}} \circ ( P_{t_{i-1},t_{i}}  - Q_{t_{i-1},t_{i}}  ) \circ P_{t_{i},t_n} h(x)\right|\nonumber\\
&+ \sup_{h\in{\rm Lip}(1)}\left| Q_{0,t_{n-1}} \circ ( P_{t_{n-1},t_{n}}  - Q_{t_{n-1},t_{n}}  )h(x)\right|.
\end{align}

We first bound the last term in \eqref{framework}. Since $h\in {\rm Lip}(1)$, one can derive from \eqref{3} and
\eqref{e:Y-moment1} with $\kappa=1$ that
\begin{align}\label{secondterm}
&\sup_{h\in {\rm Lip}(1)} \left|Q_{0,t_{n-1}} \circ  ( P_{t_{n-1},t_{n}}  - Q_{t_{n-1},t_{n}}  )h(x)\right|\nonumber\\
=&\sup_{h\in {\rm Lip}(1)}\left|\mathbb{E}[h(X_{t_{n-1},t_{n}}^{Y_{t_{n-1}}^{x}})]-\mathbb{E}[h(Y_{t_{n-1},t_{n}}^{Y_{t_{n-1}}^{x}})]\right|\nonumber\\
\leq&\mathbb{E}\left|X_{t_{n-1},t_{n}}^{Y_{t_{n-1}}^{x}}-Y_{t_{n-1},t_{n}}^{Y_{t_{n-1}}^{x}}\right|
\leq C d^{\frac{5}{2} } (1+\mathbb{E}|Y_{t_{n-1}}^{x}|)\gamma_{n}^{1+\frac{1}{\alpha}} \nonumber  \\
\leq& Cd^{5}(1+|x|)\gamma_{n}^{1+\frac{1}{\alpha}}.
\end{align}

For the first sum in \eqref{framework}, we use Lemma \ref{pro:erg-X} and the Markov property to see
that  there exists a  constant $C>0$ independent of $d$ such that for any $x,y\in \R^d$,
\begin{eqnarray}\label{ergodicity1}
|P_{t_i,t_j}h(x)-P_{t_i,t_j}h(y)| &\leq& Ce^{ d^{4-2\alpha} 2^{d}}|x-y|e^{-d^{2\alpha-4}2^{-d}
e^{-2^{d} d^{4-2\alpha}}(t_{j}-t_{i})} ,
\end{eqnarray}
It follows from this, \eqref{3}, and \eqref{e:Y-moment1} that
\begin{eqnarray*}
&& \sup_{h\in {\rm Lip}(1)} \left| Q_{0,t_{i-1}} \circ ( P_{t_{i-1},t_{i}}  - Q_{t_{i-1},t_{i}}  ) \circ P_{t_i,t_n} h(x) \right|  \\
&=& \sup_{h\in {\rm Lip}(1)} \left|\mathbb{E}\left[P_{t_{i},t_{n}}h\left(X_{t_{i-1},t_{i}}^{Y_{t_{i-1}}^{x}}\right)\right]-
\mathbb{E}\left[P_{t_{i},t_{n}}h\left(Y_{t_{i-1},t_{i}}^{Y_{t_{i-1}^{x}}}\right)\right]\right| \nonumber  \\
&\leq& Ce^{ d^{4-2\alpha} 2^{d}} e^{-d^{2\alpha-4}2^{-d}
e^{-2^{d} d^{4-2\alpha}}(t_{n}-t_{i})} \mathbb{E}\left|X_{t_{i-1},t_{i}}^{Y_{t_{i-1}}^{x}}-Y_{t_{i-1},t_{i}}^{Y_{t_{i-1}^{x}}}\right|\\
&\leq&  Ce^{ d^{4-2\alpha} 2^{d}} e^{-d^{2\alpha-4}2^{-d}
e^{-2^{d} d^{4-2\alpha}}(t_{n}-t_{i})}d^{\frac{5}{2}}
\left(1+\mathbb{E}\left|Y_{t_{i-1}}^{x}\right|\right)\gamma_{i}^{1+\frac{1}{\alpha}}\\
&\leq& Ce^{ d^{4-2\alpha} 2^{d}} e^{-d^{2\alpha-4}2^{-d}
e^{-2^{d} d^{4-2\alpha}}(t_{n}-t_{i})}d^{5}
\left(1+\left|x\right|\right)\gamma_{i}^{1+\frac{1}{\alpha}}.
\end{eqnarray*}
Recall that $\rho= d^{2\alpha-4}2^{-d} e^{-2^{d} d^{4-2\alpha}},$ we have
\begin{align*}
&\sum_{i=1}^{n-1}\sup_{h\in{\rm Lip}(1)}\left|Q_{0,t_{i-1}} \circ ( P_{t_{i-1},t_{i}}  - Q_{t_{i-1},t_{i}}  ) \circ P_{t_{i},t_n} h(x)\right|\\
\leq&Ce^{d^{4-2\alpha} 2^{d}}d^{5}(1+|x|)
\sum_{i=1}^{n-1}e^{-\rho
(t_{n}-t_{i})}\gamma_{i}^{1+\frac{1}{\alpha}}
\leq\frac{C}{\rho-\omega}
e^{\frac{\rho-\omega}{2}\gamma_{1}}
e^{d^{4-2\alpha} 2^{d}}
d^{5} (1+|x|)\gamma_{n}^{\frac{1}{\alpha}},
\end{align*}
where the last inequality follows from \eqref{convergence1}.

Combining all of the above, we derive that  there exists a constant $C>0$ independent of $d$ and $n$ such that
\begin{align*}
W_1(\mathcal{L}(X_{t_n}^x), \mathcal{L}(Y_{t_n}^x))\leq&
\frac{C}{\rho-\omega}e^{\frac{\rho-\omega}{2}\gamma_{1}}
e^{d^{4-2\alpha}2^{d}}d^{5}(1+|x|)\gamma_{n}^{\frac{1}{\alpha}}.
\end{align*}
This proves \eqref{Liptheorem}.

Now we prove \eqref{LiptheoremB}. Notice that
\begin{align}\label{invarergo}
W_{1}(\nu,\mathcal{L}(X_{t_n}^x))
=&\int_{\mathbb{R}^{d}}W_{1}(\mathcal{L}(X_{t_n}^x), \mathcal{L}(X_{t_n}^y))\nu(\dif y)\nonumber\\
\leq& Ce^{d^{4-2\alpha} 2^{d}}e^{-\rho t_{n}} \int_{\mathbb{R}^{d}}|x-y|\nu(\dif y)
\leq Ce^{d^{4-2\alpha} 2^{d}} d(1+|x|)e^{-\rho t_{n}},
\end{align}
where the first inequality in (\ref{invarergo}) is from Lemma \ref{pro:erg-X}. Then, by \eqref{Liptheorem} and \eqref{convergence2},
\begin{align*}
W_{1}(\nu, \mathcal{L}(Y_{t_n}^x))\leq& W_{1}(\nu,\mathcal{L}(X_{t_n}^x))+W_{1}(\mathcal{L}(X_{t_n}^x), \mathcal{L}(Y_{t_n}^x))\\
\leq&\frac{C}{\rho-\omega}e^{\frac{\rho-\omega}{2}\gamma_1} e^{d^{4-2\alpha} 2^{d}}d^{5}(e^{-\rho t_{n}}+\gamma_{n}^{\frac{1}{\alpha}}) (1+|x|) \\
\leq&\frac{C}{\rho-\omega}e^{\frac{\rho-\omega}{2}\gamma_1} e^{d^{4-2\alpha} 2^{d}}d^{5}\gamma_{n}^{\frac{1}{\alpha}} (1+|x|).
\end{align*}
The proof is complete.
\qed

\subsection{Proof of Theorem \ref{thm:LXLY1}}
According to \eqref{framework}, we have
\begin{align}\label{frame1}
W_1(\mathcal{L}(X_{t_n}^x), \mathcal{L}(Y_{t_n}^x))
\leq& \sum_{i=1}^{n-1}\sup_{h\in{\rm Lip}(1)}\left|Q_{0,t_{i-1}} \circ ( P_{t_{i-1},t_{i}}  - Q_{t_{i-1},t_{i}}  ) \circ P_{t_{i},t_n} h(x)\right| \nonumber \\
&+ \sup_{h\in{\rm Lip}(1)}\left| Q_{0,t_{n-1}} \circ ( P_{t_{n-1},t_{n}}  - Q_{t_{n-1},t_{n}}  )h(x)\right|.
\end{align}

For the first sum on the right hand side of (\ref{frame1}), we bound it as follows:%\footnote{\textcolor{blue} {Is $\sup_{h\in{\rm Lip}(1)}$ missing?}}
\begin{eqnarray} \label{e:R-12}
&&  \sum_{i=1}^{n-1} \sup_{h\in{\rm Lip}(1)} \left| Q_{0,t_{i-1}} \circ ( P_{t_{i-1},t_{i}}- Q_{t_{i-1},t_{i}}  ) \circ P_{t_i,t_n} h(x) \right|   \nonumber  \\
&\leq&  \sum_{i=1}^{n^*} \sup_{h\in{\rm Lip}(1)} \left| Q_{0,t_{i-1}} \circ ( P_{t_{i-1},t_{i}}  - Q_{t_{i-1},t_{i}}  ) \circ P_{t_i,t_n} h(x) \right|  \nonumber  \\
&& +  \sum_{i=n^*+1}^{n-1} \sup_{h\in{\rm Lip}(1)} \left| Q_{0,t_{i-1}} \circ ( P_{t_{i-1},t_{i}}  - Q_{t_{i-1},t_{i}}  ) \circ P_{t_i,t_n} h(x) \right|   \nonumber  \\
&=:& \mathcal{R}_1 + \mathcal{R}_2,
\end{eqnarray}
recall $n^*=\max\{ i: t_n - t_i >1 \}$.

Let us first give the estimate for $\mathcal{R}_1$, in which $t_n-t_i>1$, recall that $\rho= d^{2\alpha-4}2^{-d} e^{-2^{d} d^{4-2\alpha}},$
then one can derive from \eqref{ergodicity1} that
\begin{eqnarray*}
&& \sup_{h\in {\rm Lip}(1)} | Q_{0,t_{i-1}} \circ ( P_{t_{i-1},t_{i}}  - Q_{t_{i-1},t_{i}}  ) \circ P_{t_i,t_n} h(x) |  \\
&=& \sup_{h\in {\rm Lip}(1)} | Q_{0,t_{i-1}} \circ ( P_{t_{i-1},t_{i}} - Q_{t_{i-1},t_{i}}  ) \circ P_{t_i,t_{i}+1} \circ P_{t_{i}+1,t_n} h(x) |  \nonumber  \\
&=& \sup_{\{g=P_{t_{i}+1,t_n} h: h\in {\rm Lip}(1)\}} | Q_{0,t_{i-1}} \circ ( P_{t_{i-1},t_{i}} - Q_{t_{i-1},t_{i}}  ) \circ P_{t_i,t_{i}+1} g(x) |  \nonumber  \\
&\le& \sup_{g\in {\rm Lip}(\ell)} | Q_{0,t_{i-1}} \circ ( P_{t_{i-1},t_{i}} - Q_{t_{i-1},t_{i}}  ) \circ P_{t_i,t_{i}+1} g(x) |  \nonumber  \\
&\le & Ce^{d^{4-2\alpha} 2^{d}}e^{-\rho(t_{n}-t_{i}-1)}  \sup_{ g\in {\rm Lip}(1) }  | Q_{0,t_{i-1}} \circ ( P_{t_{i-1},t_{i}} - Q_{t_{i-1},t_{i}}  ) \circ P_{t_i,t_i+1} g(x)|,
\end{eqnarray*}
where $\ell=Ce^{d^{4-2\alpha} 2^{d}}e^{-\rho(t_{n}-t_{i}-1)}$ and the first inequality follows from \eqref{ergodicity1}.
Then, it follows from Lemmas \ref{lem:Y-moment}, \ref{pro:P-Qf1} and \ref{pro:Pt} that
\begin{eqnarray*}
&& \sup_{h\in {\rm Lip}(1)} | Q_{0,t_{i-1}} \circ ( P_{t_{i-1},t_{i}}  - Q_{t_{i-1},t_{i}}  ) \circ P_{t_i,t_n} h(x) |  \\
&\leq& Ce^{d^{4-2\alpha} 2^{d}}e^{-\rho(t_{n}-t_{i}-1)}
 d^{5\vee(3\kappa)\vee(4+\frac{3}{2\kappa}) \vee (\frac{3}{2}+\frac{3}{2\kappa}+\frac{3\kappa}{2}
)}\\
 &&\qquad\qquad(1+|x|^{\kappa}) \sup_{g\in {\rm Lip}(\ell)} (\| \nabla P_{t_i,t_i+1}g\|_{\infty} + \| \nabla^2
 P_{t_i,t_i+1}g\|_{{\rm HS},\infty}) \gamma_i^{2+\frac{1}{\alpha}-\frac{1}{\kappa}}\nonumber \\
&\leq& Ce^{d^{4-2\alpha} 2^{d}}e^{-\rho(t_{n}-t_{i}-1)}
d^{5\vee(3\kappa)\vee(4+\frac{3}{2\kappa}) \vee (\frac{3}{2}+\frac{3}{2\kappa}+\frac{3\kappa}{2})}
(1+|x|^{\kappa})
\gamma_i^{2+\frac{1}{\alpha}-\frac{1}{\kappa}} ,
\end{eqnarray*}
which implies that there exists a positive constant $C$ independent of $d$ and $n$ such that
\begin{eqnarray}\label{e:R-1}
\mathcal{R}_1
&=&\sum_{i=1}^{n^*}\sup_{h\in{\rm Lip}(1)} \left|Q_{0,t_{i-1}} \circ (P_{t_{i-1},t_{i}}-Q_{t_{i-1},t_{i}}) \circ P_{t_i,t_n} h(x) \right| \nonumber  \\
&\leq& Ce^{d^{4-2\alpha} 2^{d}}
 d^{5\vee(3\kappa)\vee(4+\frac{3}{2\kappa}) \vee (\frac{3}{2}+\frac{3}{2\kappa}+\frac{3\kappa}{2})}
(1+|x|^{\kappa})\sum_{i=1}^{n^*} e^{-\rho (t_{n}-t_{i}-1)  } \gamma_i^{2+\frac{1}{\alpha}-\frac{1}{\kappa}} \nonumber \\
&\leq&\frac{ Ce^{\rho}}{\rho-\omega}e^{\frac{\rho-\omega}{2}\gamma_1} e^{d^{4-2\alpha} 2^{d}}
d^{5\vee(3\kappa)\vee(4+\frac{3}{2\kappa}) \vee (\frac{3}{2}+\frac{3}{2\kappa}+\frac{3\kappa}{2}
)}
(1+|x|^{\kappa}) \gamma_{n}^{1+\frac{1}{\alpha}-\frac{1}{\kappa}},
\end{eqnarray}
where the last inequality holds from \eqref{convergence1}.

For $\mathcal{R}_2$, it follows from Lemmas \ref{lem:Y-moment}, \ref{pro:P-Qf1}, and \ref{pro:Pt} that
\begin{eqnarray}\label{e:R-2}
\mathcal{R}_2
&=&\sum_{i=n^*+1}^{n-1} \sup_{h\in{\rm Lip}(1)}\left| Q_{0,t_{i-1}} \circ ( P_{t_{i-1},t_{i}}  - Q_{t_{i-1},t_{i}}  ) \circ P_{t_i,t_n} h(x) \right|   \nonumber  \\
&\leq&  C e^{\rho(1+\gamma_{1})}
 d^{5\vee(3\kappa)\vee(4+\frac{3}{2\kappa}) \vee (\frac{3}{2}+\frac{3}{2\kappa}+\frac{3\kappa}{2}
)} (1+|x|^{\kappa})\sum_{i=n^*+1}^{n-1}(t_n-t_i)^{-\frac{1}{\alpha}} \gamma_i^{2+\frac{1}{\alpha}-\frac{1}{\kappa}}\nonumber\\
&\leq& C e^{\rho(1+\gamma_{1})}
 d^{5\vee(3\kappa)\vee(4+\frac{3}{2\kappa}) \vee (\frac{3}{2}+\frac{3}{2\kappa}+\frac{3\kappa}{2}
)}
(1+|x|^{\kappa})
\gamma_{n}^{1+\frac{1}{\alpha}-\frac{1}{\kappa}},
\end{eqnarray}
where the constant $C$ is independent of $d$ and the last inequality is from \eqref{convergence3}. Hence, \eqref{e:R-12}, \eqref{e:R-1},
and \eqref{e:R-2} imply that
\begin{align}\label{e:ss-2}
&\sum_{i=1}^{n-1} \sup_{h\in{\rm Lip}(1)}\left| Q_{0,t_{i-1}} \circ ( P_{t_{i-1},t_{i}}  - Q_{t_{i-1},t_{i}}  ) \circ P_{t_i,t_n} h(x) \right|\nonumber\\
\leq&\frac{C}{\rho-\omega}e^{\frac{\rho-\omega}{2}\gamma_1} e^{d^{4-2\alpha} 2^{d}}
d^{5\vee(3\kappa)\vee(4+\frac{3}{2\kappa}) \vee (\frac{3}{2}+\frac{3}{2\kappa}+\frac{3\kappa}{2}
)}
(1+|x|^{\kappa}) \gamma_{n}^{1+\frac{1}{\alpha}-\frac{1}{\kappa}}.
\end{align}

In addition, the second term in (\ref{frame1}) is directly obtained from \eqref{secondterm}, that is,
\begin{align*}
\sup_{h\in{\rm Lip}(1)}\left|Q_{0,t_{n-1}}\circ ( P_{t_{n-1},t_{n}}  - Q_{t_{n-1},t_{n}}  )h(x)\right|
\leq& Cd^{5} (1+|x|)\gamma_{n}^{1+\frac{1}{\alpha}}.
\end{align*}
Therefore, there exists a positive constant $C$ independent of $d$ such that
\begin{eqnarray*}
W_{1}(\mathcal{L}(X_{t_n}^x), \mathcal{L}(Y_{t_n}^x))
\leq&\frac{C}{\rho-\omega}e^{\frac{\rho-\omega}{2}\gamma_1}e^{d^{4-2\alpha} 2^{d}}
d^{5\vee(3\kappa)\vee(4+\frac{3}{2\kappa}) \vee (\frac{3}{2}+\frac{3}{2\kappa}+\frac{3\kappa}{2}
)} (1+|x|^{\kappa}) \gamma_{n}^{1+\frac{1}{\alpha}-\frac{1}{\kappa}}.
\end{eqnarray*}
This proves \eqref{e:w-distance1}.

Now we prove (\ref{e:w-distance1B}). By \eqref{invarergo},  \eqref{e:w-distance1} and \eqref{convergence2}, we have
\begin{align*}
W_{1}(\nu, \mathcal{L}(Y_{t_n}^x))\leq& W_{1}(\nu,\mathcal{L}(X_{t_n}^x))+W_{1}(\mathcal{L}(X_{t_n}^x), \mathcal{L}(Y_{t_n}^x))\\
\leq&\frac{C}{\rho-\omega}e^{\frac{\rho-\omega}{2}\gamma_1}e^{d^{4-2\alpha} 2^{d}}
 d^{5\vee(3\kappa)\vee(4+\frac{3}{2\kappa}) \vee (\frac{3}{2}+\frac{3}{2\kappa}+\frac{3\kappa}{2}
)}
(e^{-\rho t_{n}}+\gamma_{n}^{1+\frac{1}{\alpha}-\frac{1}{\kappa}}) (1+|x|^{\kappa}) \\
\leq&\frac{C}{\rho-\omega}e^{\frac{\rho-\omega}{2}\gamma_1}e^{d^{4-2\alpha} 2^{d}}
d^{5\vee(3\kappa)\vee(4+\frac{3}{2\kappa}) \vee (\frac{3}{2}+\frac{3}{2\kappa}+\frac{3\kappa}{2}
)} \gamma_{n}^{1+\frac{1}{\alpha}-\frac{1}{\kappa}} (1+|x|^{\kappa}) .
\end{align*}
The proof is complete.
\qed

\subsection{Proof of Theorem \ref{thm:LXLY}}
 By \eqref{frame1}, we can write
\begin{align*}
W_1(\mathcal{L}(X_{t_n}^x), \mathcal{L}(\tilde{Y}_{t_n}^x))
\leq&
\sum_{i=1}^{n-1}\sup_{h\in{\rm Lip}(1)}\left|\tilde{Q}_{0,t_{i-1}} \circ ( P_{t_{i-1},t_{i}}  - \tilde{Q}_{t_{i-1},t_{i}}  ) \circ P_{t_{i},t_n} h(x)\right|\\
&+\sup_{h\in{\rm Lip}(1)}\left|\tilde{Q}_{0,t_{n-1}} \circ ( P_{t_{n-1},t_{n}}  - \tilde{Q}_{t_{n-1},t_{n}}  )h(x)\right|,
\end{align*}
whereas by the same argument as the proof of \eqref{e:ss-2} with Lemma \ref{pro:P-Qf1} replaced by Lemma \ref{pro:P-Qf}, we have
\begin{align*}
&\sum_{i=1}^{n-1} \sup_{h\in{\rm Lip}(1)} \left|\tilde{Q}_{0,t_{i-1}} \circ ( P_{t_{i-1},t_{i}}  - \tilde{Q}_{t_{i-1},t_{i}}  ) \circ P_{t_i,t_n} h(x) \right|\\
\leq&\frac{C}{\rho-\omega}e^{\frac{\rho-\omega}{2}\gamma_1}e^{d^{4-2\alpha} 2^{d}}
d^{\frac{5(1+\kappa)}{2\kappa}
\vee(\frac{11-\alpha}{2}
)}
(1+|x|) \gamma_{n}^{\frac{2-\alpha}{\alpha}}.
\end{align*}
In addition,  \eqref{e:Y-moment} and \eqref{2} give
\begin{align*}
& \sup_{h\in{\rm Lip}(1)}\left|\tilde{Q}_{0,t_{n-1}} \circ ( P_{t_{n-1},t_{n}}  - \tilde{Q}_{t_{n-1},t_{n}}  )h(x)\right|\\
\leq& \sup_{h\in{\rm Lip}(1)}\left|\mathbb{E}
\left[h\left(X_{t_{n-1},t_{n}}^{\tilde{Y}_{t_{n-1}}^{x}}\right)
-h\left(\tilde{Y}_{t_{n-1}}^{x}\right)\right]\right|
+ \sup_{h\in{\rm Lip}(1)} \left|\mathbb{E}\left[h\left(\tilde{Y}_{t_{n-1},t_{n}}^{\tilde{Y}_{t_{n-1}}^{x}}\right)-h\left(\tilde{Y}_{t_{n-1}}^{x}\right)\right]\right|\\
\leq&\mathbb{E}\left|X_{t_{n-1},t_{n}}^{\tilde{Y}_{t_{n-1}}^{x}}-\tilde{Y}_{t_{n-1}}^{x}\right|
+\mathbb{E}\left|\gamma_{n}b(\tilde{Y}_{t_{n-1}})+\frac{\gamma_{n}^{\frac{1}{\alpha}}}{\beta}A\tilde{Z}_{n}\right|\\
\leq&Cd^{\frac{5}{2}} \left(1+\mathbb{E}|\tilde{Y}_{t_{n-1}}^{x}|\right)\gamma_{n}^{\frac{1}{\alpha}}\leq Cd^{5} (1+|x|)\gamma_{n}^{\frac{1}{\alpha}}.
\end{align*}
Combining the above, we obtain
\begin{eqnarray*}
W_{1}(\mathcal{L}(X_{t_n}^x), \mathcal{L}(\tilde{Y}_{t_n}^x))
\leq&\frac{C}{\rho-\omega}e^{\frac{\rho-\omega}{2}\gamma_1}e^{d^{4-2\alpha} 2^{d}}
d^{5}
(1+|x|) \gamma_{n}^{\frac{2-\alpha}{\alpha}}.
\end{eqnarray*}
This proves \eqref{e:w-distance}.

Next we prove \eqref{e:w-distanceB}. By \eqref{invarergo}, \eqref{e:w-distance} and \eqref{convergence2}, we have
\begin{align*}
W_{1}(\nu, \mathcal{L}(\tilde{Y}_{t_n}^x))\leq& W_{1}(\nu,\mathcal{L}(X_{t_n}^x))+W_{1}(\mathcal{L}(X_{t_n}^x), \mathcal{L}(\tilde{Y}_{t_n}^x))\\
\leq&\frac{C}{\rho-\omega}e^{\frac{\rho-\omega}{2}\gamma_1}e^{d^{4-2\alpha} 2^{d}}d^{5}
 (e^{-\rho t_{n}}+\gamma_{n}^{\frac{2-\alpha}{\alpha}}) (1+|x|)\\
\leq&\frac{C}{\rho-\omega}e^{\frac{\rho-\omega}{2}\gamma_1}e^{d^{4-2\alpha} 2^{d}}d^{5}\gamma_{n}^{\frac{2-\alpha}{\alpha}} (1+|x|).
\end{align*}
The proof is complete.
\qed

\section{Proof of The Auxiliary Lemmas in Section \ref{sec:proof-thm}} \label{APP:ss}

\subsection{Proof of Lemma \ref{lem:Y-moment}}

For any $\kappa \in [1,\alpha)$, denote
$$V_{\kappa}(x)=(1+|x|^2)^{\frac{\kappa}{2}},$$
then $|x|^{\kappa}\leq V_{\kappa}(x)\leq1+|x|^{\kappa}$.

We first prove \eqref{X-moment}. Recall that the infinitesimal generator $\mathcal{A}$ of the process $(X_{t})_{t\ge 0}$ is
\begin{eqnarray}\label{operator}
\mathcal{L} f(x)
&=& \langle b(x), \nabla f(x) \rangle + \Delta^{\frac{\alpha}{2},A}f(x),
\end{eqnarray}
where
\begin{align}\label{frac}
\Delta^{\frac{\alpha}{2},A}f(x)=\int_{ \R^d \setminus \{0\} } [f(x+Az)-f(x)-\langle Az, \nabla f(x) \rangle 1_{\{|z|\leq 1\}}] \frac{d_{\alpha}}{|z|^{d+\alpha}} \dif z.
\end{align}
By It\^{o}'s formula and the same argument as the proof of \cite[A.3]{Chen2022stableEM}, we have
\begin{align*}
\frac{\dif \mathbb{E}\left[V_{\kappa}\left(X_{t}^{x}\right)\right]}{\dif t}=&\mathbb{E}\left[\left\langle b\left(X_{t}^{x}\right),
\nabla V_{\kappa}\left(X_{t}^{x}\right)\right\rangle\right]+\mathbb{E}\left[\Delta^{\frac{\alpha}{2}}V_{\kappa}\left(X_{t}^{x}\right)\right]\\
\leq&-\frac{\theta_{1}}{2}\mathbb{E}\left[V_{\kappa}\left(X_{t}^{x}\right)\right]+C_{1},
\end{align*}
where
\begin{align*}
C_{1}=&\theta_{1}\kappa+\kappa K+\theta_{1}^{1-\kappa}|b(0)|^{\kappa}\\
&+\frac{d_{\alpha}\kappa(3-\kappa)\sqrt{d}\sigma_{d-1}}{2(2-\alpha)}\|A\|_{{\rm op}}^{2}+\frac{d_{\alpha}\kappa\sigma_{d-1}}
{\alpha-\kappa}\|A\|_{{\rm op}}^{\kappa}
+\left(\frac{\theta_{1}}{4}\right)^{1-\kappa}\left(\frac{d_{\alpha}\sigma_{d-1}}{\alpha-1}\|A\|_{{\rm op}}\right)^{\kappa}.
\end{align*}
This inequality, together with $X_{0}=x$, implies
\begin{align*}
\mathbb{E}\left[V_{\kappa}\left(X_{t}^{x}\right)\right]\leq e^{-\frac{\theta_{1}}{2}t}V_{\kappa}(x)+C_{1}\int_{0}^{t}e^{-\frac{\theta_{1}}{2}(t-s)}\dif s
\leq V_{\kappa}(x)+\frac{2C_{1}}{\theta_{1}}.
\end{align*}
Hence, the desired result follows from the fact
\begin{align*}
C_{1}\leq Cd^{\frac{5\vee(3\kappa)}{2}}
\end{align*}
for some constant $C>0$, which is independent of $d$.

For the remaining two inequalities, we only give the proof of \eqref{e:Y-moment1}. Inequality \eqref{e:Y-moment} can be proved
in the same way.

By the same argument as the proof of Lemma 1.8 in \cite{Chen2022stableEM}, we have
\begin{align*}
\mathbb{E}\left[V_{\kappa}\left(Y_{t_{n+1}}^{x}\right)\right]
\leq\left(1-\frac{\theta_{1}}{2}\gamma_{n+1}\right)
\mathbb{E}[V_{\kappa}(Y_{t_{n}}^x)]+C_{2}\gamma_{n+1},
\end{align*}
where
\begin{align*}
C_{2}=&\frac{\theta_{1}\kappa}{2} \left(\frac{2\gamma_{n+1}|b(0)|^{2}}{\theta_{1}}
+2\gamma_{n+1}^{2}|b(0)|^{2}+1+2\gamma_{n+1} K\right)+ \frac{\kappa |b(0)|^{2}}{\theta_{1}}+2\kappa \gamma_{n+1} |b(0)|^{2}\\
&+\kappa K+d_{\alpha}\kappa\sigma_{d-1}\left(\frac{(3-\kappa)\sqrt{d}}{2(2-\alpha)}\|A\|_{{\rm op}}^{2}
+\frac{1}{\alpha-\kappa}\|A\|_{{\rm op}}^{\kappa}+\frac{|b(0)|^{\kappa-1}}{\alpha-1}\|A\|_{{\rm op}}\right)\\
&\qquad+d_{\alpha}\kappa\sigma_{d-1}\frac{\mathbb{E}|Z_{1}|^{\kappa-1}}{\alpha-1}\|A\|_{{\rm op}}+\left(\frac{d_{\alpha}\sigma_{d-1}
(1+ L^{\kappa-1})}{\alpha-1}\|A\|_{{\rm op}}\right)^{\kappa}
\left(\frac{2}{\theta_{1}}\right)^{\kappa-1}.
\end{align*}
By induction, we further have
\begin{align*}
\mathbb{E}\left[V_{\kappa}\left(Y_{t_{n+1}}^{x}\right)\right]
\leq&\prod_{i=1}^{n+1}\left(1-\frac{\theta_{1}}{2}
\gamma_{i}\right)V_{\kappa}(x)
+C_{2}\sum_{i=1}^{n+1}\prod_{j=i+1}^{n+1}
\left(1-\frac{\theta_{1}}{2}\gamma_{j}\right)\gamma_{i}\\
\leq&V_{\kappa}(x)+C_{2}\sum_{i=1}^{n+1}\prod_{j=i+1}^{n+1}
\left(1-\frac{\theta_{1}}{2}\gamma_{j}\right)\gamma_{i},
\end{align*}
here we have used the convention that $\prod_{j=n+2}^{n+1}\left(1-\frac{\theta_{1}}{2}\gamma_{j}\right)=1$. In addition,
\begin{align*}
\frac{\theta_{1}}{2}\sum_{i=1}^{n+1}\prod_{j=i+1}^{n+1}\left(1-\frac{\theta_{1}}{2}\gamma_{j}\right)\gamma_{i}
=1-\prod_{j=1}^{n+1}\left(1-\frac{\theta_{1}}{2}\gamma_{j}\right)\leq 1,
\end{align*}
which implies
\begin{align*}
\mathbb{E}\left[V_{\kappa}\left(Y_{t_{n+1}}^{x}\right)\right]
\leq&V_{\kappa}(x)+\frac{2C_{2}}{\theta_{1}},
\end{align*}
and $V_{\kappa}(x)\leq1+|x|^{\kappa}$, the desired result follows from the fact
\begin{align*}
C_{2}\leq Cd^{\frac{5\vee(3\kappa)}{2}}
\end{align*}
for some constant $C>0$, which is independent of $d$.
\qed

\subsection{Proof of Lemma \ref{lem:mom}}
The first inequality follows immediately from
\begin{align*}
    \mathbb{E} |Y_{t_{n-1},t_{n}}^{x}-x|^{\kappa}
    =& \mathbb{E}\left|\gamma_{n} b(x) + A\left(Z_{t_{n}}-Z_{t_{n-1}}\right)\right|^{\kappa}\\
    \leq& 2\left[\gamma_{n}^{\kappa}|b(x)|^{\kappa}+\|A\|_{{\rm op}}^{\kappa} \mathbb{E}|Z_{\gamma_{n}}|^{\kappa}\right]
    \leq Cd^{\kappa}(1+|x|^{\kappa})\gamma_{n}^{\kappa/\alpha}.
\end{align*}

By the H\"{o}lder inequality and \eqref{X-moment}, we obtain
\begin{align*}
    \mathbb{E}|X_{s,t}^{x}-x|^{\kappa}
    &=\mathbb{E}\left|\int_{s}^{t}b\left(X_{s,r}^{x}\right)\,\dif r+A\left(Z_{t}-Z_{s}\right)\right|^{\kappa}\\
    &\leq 2\mathbb{E}\left|\int_{s}^{t}b\left(X_{s,r}^{x}\right)\,\dif r\right|^{\kappa}+2\|A\|_{{\rm op}}^{\kappa}|Z_{t-s}|^{\kappa}\\
    &\leq 2(t-s)^{\kappa-1}\int_{s}^{t}\mathbb{E}\left|b\left(X_{s,r}^{x}\right)\right|^{\kappa}\,\dif r+2d^{\frac{\kappa}{2}}
    \mathbb{E}|Z_{1}|^{\kappa}(t-s)^{\frac{\kappa}{\alpha}}\\
    &\leq Cd^{\frac{5\vee(3\kappa)}{2}}(1+|x|^{\kappa})(t-s)^{\kappa/\alpha},
\end{align*}
which implies the second inequality.

For the last inequality, notice that $\gamma_{n}=t_{n}-t_{n-1}$, the H\"{o}lder inequality, \eqref{Linear} and \eqref{2} imply
\begin{align*}
    \mathbb{E} |X_{t_{n-1},t_{n}}^{x}-Y_{t_{n-1},t_{n}}^{x}|^{\kappa}
    &=\mathbb{E} \left|\int_{t_{n-1}}^{t_{n}}\left[b\left(X_{t_{n-1},s}^{x}\right)-b(x)\right]\,\dif s\right|^{\kappa}\\
    &\leq \gamma_{n}^{\kappa-1}\int_{t_{n-1}}^{t_{n}}\mathbb{E}\left|b\left(X_{t_{n-1},s}^{x}\right)-b(x)\right|^{\kappa}\,\dif s\\
    &\leq L\gamma_{n}^{\kappa-1}\int_{t_{n-1}}^{t_{n}}\mathbb{E}\left|X_{t_{n-1},s}^{x}-x\right|^{\kappa}\,\dif s\\
    &\leq Cd^{\frac{5\vee(3\kappa)}{2}}(1+|x|^{\kappa})\gamma_{n}^{\kappa-1}\int_{t_{n-1}}^{t_{n}}(s-t_{n-1})^{\frac{\kappa}{\alpha}}\,\dif s\\
    &=Cd^{\frac{5\vee(3\kappa)}{2}}(1+|x|^{\kappa})\gamma_{n}^{\kappa+\frac{\kappa}{\alpha}}.
\end{align*}
\qed

\subsection{Proof of Lemma \ref{pro:P-Qf1}}

By the Taylor expansion formula, we have
\begin{align*}
    &\mathbb{E} f\left(X_{t_{k-1},t_{k}}^{x}\right)-\mathbb{E} f\left(Y_{t_{k-1},t_{k}}^{x}\right)\\
    &=\mathbb{E}\left\langle\nabla f\left(Y_{t_{k-1},t_{k}}^{x}\right),X_{t_{k-1},t_{k}}^{x}-Y_{t_{k-1},t_{k}}^{x}\right\rangle\\
    &+\mathbb{E}\int_{0}^{1}\!\!\left\langle\nabla f\left(Y_{t_{k-1},t_{k}}^{x}+r\left(X_{t_{k-1},t_{k}}^{x}-Y_{t_{k-1},t_{k}}^{x}\right)\right)\!-\!\nabla f\left(Y_{t_{k-1},t_{k}}^{x}\right),X_{t_{k-1},t_{k}}^{x}-Y_{t_{k-1},t_{k}}^{x}\right\rangle\dif r\\
    &=:\mathcal{I}+\mathcal{II}+\mathcal{III},
\end{align*}
where
\begin{align*}
\mathcal{I}=&\mathbb{E}\left\langle\nabla f\left(x+\gamma_{k}b(x)+A(Z_{t_{k}}-Z_{t_{k-1}})\right)-\nabla f\left(x+\gamma_{k}b(x)\right),X_{t_{k-1},t_{k}}^{x}-Y_{t_{k-1},t_{k}}^{x}\right\rangle, \\
\mathcal{II}=&\mathbb{E}\left\langle\nabla f\left(x+\gamma_{k}b(x)\right),X_{t_{k-1},t_{k}}^{x}-Y_{t_{k-1},t_{k}}^{x}\right\rangle, \\
\mathcal{III}=&\mathbb{E}\int_{0}^{1}
\left\langle\nabla f\left(Y_{t_{k-1},t_{k}}^{x}+r\left(X_{t_{k-1},t_{k}}^{x}-Y_{t_{k-1},t_{k}}^{x}\right)\right)\right.\\
& \left.\qquad \qquad \qquad -\nabla f\left(Y_{t_{k-1},t_{k}}^{x}\right),X_{t_{k-1},t_{k}}^{x}-Y_{t_{k-1},t_{k}}^{x}\right\rangle\dif r.
\end{align*}
For the first term $\mathcal{I}$, we divide it into two parts as follows:
\begin{align*}
    \mathcal{I}
    =&\mathbb{E}\left[\left\langle\nabla f\left(x+\gamma_{k}b(x)+A(Z_{t_{k}}-Z_{t_{k-1}})\right)-\nabla f\left(x+\gamma_{k}b(x)\right),X_{t_{k-1},t_{k}}^{x}-Y_{t_{k-1},t_{k}}^{x}\right\rangle\right.\\
    &\left.\qquad\qquad\qquad\qquad\cdot\left({\bf 1}_{(0,1]}(|Z_{t_{k}}-Z_{t_{k-1}}|)+{\bf 1}_{(1,\infty)}(|Z_{t_{k}}-Z_{t_{k-1}}|)\right)\right]\\
    =&:\mathcal{I}_{1}+\mathcal{I}_{2}.
\end{align*}
According to \cite[Lemma 18]{CNXY23}, it is straightforward to derive the following estimates:
\begin{align*}
\mathds{P}\left(|Z_{\gamma_{k}}| > 1\right) \leq& C\pi^{-\frac{d}{2}}\Gamma\left(\frac{d+\alpha}{2}\right)\int_{|z|\geq\gamma_{k}^{-\frac{1}{\alpha}}}\frac{1}{|z|^{\alpha+d}}\dif z\\
\leq&C\pi^{-\frac{d}{2}}\Gamma\left(\frac{d+\alpha}{2}\right)\sigma_{d-1}\int_{\gamma_{k}^{-\frac{1}{\alpha}}}^{\infty}\frac{1}{r^{\alpha+1}}\dif r\leq C\frac{\Gamma\left(\frac{d+\alpha}{2}\right)}{\Gamma\left(\frac{d}{2}\right)}\gamma_{k}\leq Cd\gamma_{k},
\end{align*}
and for any $\lambda>\alpha$,
\begin{align*}
\mathbb{E}\left[|Z_{\gamma_{k}}|^\lambda {\bf 1}_{(0,1]}(|Z_{\gamma_{k}}|)\right]\leq& C\pi^{-\frac{d}{2}}\Gamma\left(\frac{d+\alpha}{2}\right)
\gamma_{k}^{\frac{\lambda}{\alpha}}
\int_{|z|\leq\gamma_{k}^{-\frac{1}{\alpha}}}\frac{|z|^{\lambda}}{|z|^{\alpha+d}}\dif z\\
\leq& C\pi^{-\frac{d}{2}}\Gamma\left(\frac{d+\alpha}{2}\right)\sigma_{d-1}\gamma_{k}^{\frac{\lambda}{\alpha}}\int_{0}^{\gamma_{k}^{-\frac{1}{\alpha}}}\frac{r^{\lambda}}{r^{\alpha+1}}\dif r\\
\leq& C\frac{\Gamma\left(\frac{d+\alpha}{2}\right)}{\Gamma\left(\frac{d}{2}\right)}\gamma_{k}\leq Cd\gamma_{k}.
\end{align*}
 Since $\frac{\kappa}{\kappa-1} > \alpha$, we can use the H\"{o}lder inequality and \eqref{3} to get
\begin{align*}
    |\mathcal{I}_{1}|
    &\leq \|\nabla^{2}f\|_{\mathrm{HS},\infty}\|A\|_{{\rm op}}\mathbb{E}\left[|Z_{t_{k}}-Z_{t_{k-1}}|{\bf 1}_{(0,1]}(|Z_{t_{k}}-Z_{t_{k-1}}|)\left|X_{t_{k-1},t_{k}}^{x}-Y_{t_{k-1},t_{k}}^{x}\right|\right]\\
    &\leq\|\nabla^{2}f\|_{\mathrm{HS},\infty}\sqrt{d} \left(\mathbb{E}\left[|Z_{\gamma_{k}}|^{\frac{\kappa}{\kappa-1}}{\bf 1}_{(0,1]}(|Z_{\gamma_{k}}|)\right]\right)^{\frac{\kappa-1}{\kappa}}
    \left(\mathbb{E}\left|X_{t_{k-1},t_{k}}^{x}-Y_{t_{k-1},t_{k}}^{x}\right|^{\kappa}\right)^{\frac{1}{\kappa}}\\
    &\leq C (1+|x|)\|\nabla^{2}f\|_{\mathrm{HS},\infty}
    d^{(\frac{3}{2}+\frac{3}{2\kappa})\vee (3-\frac{1}{\kappa})} \gamma_{k}^{2+\frac{1}{\alpha}-\frac{1}{\kappa}},
\end{align*}
whereas by the H\"{o}lder inequality
\begin{align*}
    |\mathcal{I}_{2}|
    &\leq2\|\nabla f\|_{\infty} \mathbb{E}\left[{\bf 1}_{(1,\infty)}(|Z_{t_{k}}-Z_{t_{k-1}}|)\left|X_{t_{k-1},t_{k}}^{x}-Y_{t_{k-1},t_{k}}^{x}\right|\right]\\
    &\leq2\|\nabla f\|_{\infty} \left(\mathbb{E}{\bf 1}_{(1,\infty)}(|Z_{\gamma_{k}}|)\right)^{\frac{\kappa-1}{\kappa}} \left(\mathbb{E}\left|X_{t_{k-1},t_{k}}^{x}-Y_{t_{k-1},t_{k}}^{x}\right|^{\kappa}\right)^{\frac{1}{\kappa}}\\
    &\leq C (1+|x|) \|\nabla f\|_{\infty}d^{(1+\frac{3}{2\kappa})\vee (\frac{5}{2}-\frac{1}{\kappa})}\gamma_{k}^{\frac{\kappa-1}{\kappa}}\gamma_{k}^{1+\frac{1}{\alpha}}.
\end{align*}
Hence, we have
\begin{align*}
    |\mathcal{I}|
    \leq C (1+|x|)\left(\|\nabla f\|_{\infty}+\|\nabla^{2}f\|_{\mathrm{HS},\infty}\right)
    d^{(\frac{3}{2}+\frac{3}{2\kappa})\vee (3-\frac{1}{\kappa})} \gamma_{k}^{2+\frac{1}{\alpha}-\frac{1}{\kappa}}.
\end{align*}

For $\mathcal{II}$, one can derive from \eqref{e:SDE}, \eqref{e:EM}, It\^{o}'s formula and \eqref{operator} that
\begin{align*}
    |\mathcal{II}|
    &=\left|\mathbb{E}\left\langle\nabla f\left(x+\gamma_{k} b(x)\right),\int_{t_{k-1}}^{t_{k}}\left[b(X_{t_{k-1},s}^{x})-b(x)\right] \dif s\right\rangle\right|\\
    &=\left|\left\langle\nabla f\left(x+\gamma_{k} b(x)\right),\int_{t_{k-1}}^{t_{k}}\mathbb{E}\left[b(X_{t_{k-1},s}^{x})-b(x)\right] \dif s\right\rangle\right|\\
    &=\left|\left\langle\nabla f\left(x+\gamma_{k} b(x)\right),\int_{t_{k-1}}^{t_{k}}\int_{t_{k-1}}^{s}\mathbb{E}\left[\mathcal{L}b(X_{t_{k-1},r}^{x})\right] \dif r\,\dif s\right\rangle\right|\\
    &\leq C\|\nabla f\|_{\infty} (1+|x|)d^{\frac{5\vee(3\kappa)}{2}} \int_{t_{k-1}}^{t_{k}}(s-t_{k-1})\dif s \nonumber \\
    &\leq C\|\nabla f\|_{\infty} (1+|x|)d^{\frac{5\vee(3\kappa)}{2}}\gamma_{k}^{2},
\end{align*}
where the last two inequality is from \eqref{Linear}, the moment estimate \eqref{X-moment}, and the inequality
\begin{align*}
\left|\mathbb{E}\left[\mathcal{L}b(X_{t_{k-1},r}^{x})\right]\right|
\leq&\mathbb{E}\left[\left|\nabla b(X_{t_{k-1},r}^{x})\right|\left|b(X_{t_{k-1},r}^{x})\right|\right]+\mathbb{E}\left|\Delta^{\frac{\alpha}{2}}b(X_{t_{k-1},r}^{x})\right|\\
\leq&C\left[1+\mathbb{E}\left|X_{t_{k-1},r}^{x}\right|+d_{\alpha}\sigma_{d-1}(\sqrt{d}\|A\|_{{\rm op}}^{2}+\|A\|_{{\rm op}})\right]\\
\leq& C(1+|x|)d^{\frac{5\vee(3\kappa)}{2}},
\end{align*}
which is proved by the same argument as in the proof of \cite[(A.2)]{Chen2022stableEM}.

Finally, by \cite[Lemma 2.6]{Chen2022stableEM} and \eqref{3}, we have
\begin{align*}
    \mathcal{III}
    &\leq C\left(\|\nabla f\|_{\infty}+\|\nabla^{2}f\|_{\mathrm{HS},\infty}\right)\mathbb{E}\left|X_{t_{k-1},t_{k}}^{x}-Y_{t_{k-1},t_{k}}^{x}\right|^{\kappa}\\
    &\leq C (1+|x|^{\kappa})  \left(\|\nabla f\|_{\infty}+\|\nabla^{2}f\|_{\mathrm{HS},\infty}\right)
    d^{\frac{5\vee(3\kappa)}{2}} \gamma_{k}^{\kappa+\frac{\kappa}{\alpha}}.
\end{align*}
This finishes the proof.
\qed
~\\

\subsection{Proof of Lemma \ref{pro:P-Qf}}

Before proving Lemma \ref{pro:P-Qf}, we first give the following lemma:

\begin{lemma}\label{laplace}
Let $\alpha\in(1,2)$ and $\Delta^{\alpha/2,A}$ be defined by \eqref{frac}. Then for any $f:\mathbb{R}^{d}\to \mathbb{R}$
satisfying $\|\nabla f\|_{\infty}<\infty$ and $\|\nabla^{2}f\|_{\mathrm{HS},\infty}<\infty$, we have
\begin{equation} \label{crucialbis}
\left|\Delta^{\alpha/2,A}f(x) -\Delta^{\alpha/2,A}f(y)\right| \leq
        \frac{d_{\alpha}\|\nabla^{2}f\|_{\mathrm{HS},\infty}\|A\|_{{\rm op}}^{2}\sigma_{d-1}}{(2-\alpha)(\alpha-1)}|x-y|^{2-\alpha},
        \quad \forall x,y\in\mathbb{R}^{d}.
\end{equation}
\end{lemma}

\begin{proof}
According to \eqref{frac} and the symmetry of the Lebesgue measure, for any $R>0$, we have
\begin{align*}
    \Delta^{\alpha/2,A}f(x)
    & =d_{\alpha}\int_{\mathds{S}^{d-1}}\int_{0}^{\infty}\frac{f(x+Ar\theta)-f(x)-Ar\left\langle\theta,\nabla f(x)\right\rangle{\bf 1}_{(0,R)}(r)}{r^{\alpha+1}}\,\dif r\,\dif\theta\\
    &= d_{\alpha}\int_{\mathds{S}^{d-1}} \int_{0}^{R}\int_{0}^{r}\frac{\left\langle A\theta,\nabla f(x+A\theta s)-\nabla f(x)\right\rangle}{r^{\alpha+1}}\,\dif s\,\dif r\,\dif \theta\\
    &\quad\mbox{}+d_{\alpha}\int_{\mathds{S}^{d-1}} \int_{R}^{\infty}\int_{0}^{r}\frac{\left\langle A\theta,\nabla f(x+A\theta s)\right\rangle}{r^{\alpha+1}}\,\dif s\,\dif r\,\dif \theta.
\end{align*}
Then, for all $x,y\in\mathbb{R}^{d}$,
\begin{align}\label{Eq:star}
&\qquad \quad \big|\Delta^{\alpha/2,A}f(x) -\Delta^{\alpha/2,A}f(y)\big| \nonumber\\
    &\qquad \ \ \ \leq d_{\alpha}\|A\|_{{\rm op}}\int_{\mathds{S}^{d-1}} \int_{0}^{R}\int_{0}^{r}\frac{\left|\nabla f(x+A\theta s)-\nabla f(x)-\nabla f(y+A\theta s)
    +\nabla f(y)\right|}{r^{\alpha+1}}\,\dif s\,\dif r\,\dif \theta \nonumber\\
    & \qquad \qquad \mbox{}+d_{\alpha}\|A\|_{{\rm op}}\int_{\mathds{S}^{d-1}} \int_{R}^{\infty}\int_{0}^{r}\frac{\left|\nabla f(x+A\theta s)-\nabla f(y+A\theta s)\right|}
    {r^{\alpha+1}}\,\dif s\,\dif r\,\dif \theta.
\end{align}
For the first integral in (\ref{Eq:star}), we have
\begin{align*}
    &\int_{\mathds{S}^{d-1}} \int_{0}^{R}\int_{0}^{r}\frac{\left|\nabla f(x+A\theta s)-\nabla f(x)-\nabla f(y+A\theta s)+\nabla f(y)\right|}{r^{\alpha+1}}\,\dif s\,\dif r\,\dif \theta\\
    &\quad\leq\int_{\mathds{S}^{d-1}} \int_{0}^{R}\int_{0}^{r}\frac{\left|\nabla f(x+A\theta s)-\nabla f(x)\right|+\left|\nabla f(y+A\theta s)-\nabla f(y)\right|}{r^{\alpha+1}}\,\dif s\,\dif r\,\dif \theta\\
    &\quad\leq 2\|\nabla^{2}f\|_{\mathrm{HS},\infty}\|A\|_{{\rm op}}\int_{\mathds{S}^{d-1}} \int_0^R \int_{0}^{r}
    \frac{s}{r^{\alpha+1}}\,\dif s\,\dif r\,\dif\theta
    =\frac{\|\nabla^{2}f\|_{\mathrm{HS},\infty}\|A\|_{{\rm op}}\sigma_{d-1}}{2-\alpha}\,R^{2-\alpha}.
\end{align*}
For the second integral in \eqref{Eq:star}, we have
\begin{align*}
    &\int_{\mathds{S}^{d-1}} \int_{R}^{\infty}\int_{0}^{r}\frac{\left|\nabla f(x+A\theta s)-\nabla f(y+A\theta s)\right|}{r^{\alpha+1}}\,\dif s\,\dif r\,\dif \theta\\
    &\quad\leq \|\nabla^{2}f\|_{\mathrm{HS},\infty}\|A\|_{{\rm op}}\int_{\mathds{S}^{d-1}}
    \int_{R}^{\infty}\int_{0}^{r}\frac{|x-y|}
    {r^{\alpha+1}}\,\dif s\,\dif r\,\dif\theta\\
    &\quad= \frac{\|\nabla^{2}f\|_{\mathrm{HS},\infty}\|A\|_{{\rm op}}\sigma_{d-1}}{\alpha-1}\,|x-y|R^{1-\alpha}.
\end{align*}
Hence, the assertion follows upon taking $R=|x-y|$.
\end{proof}

Now, we proceed to prove Lemma \ref{pro:P-Qf}.

\begin{proof}[Proof of Lemma \ref{pro:P-Qf}]
From \eqref{e:SDE} and \eqref{e:EM2}, we see
\begin{align*}
    &\mathbb{E} [f(X_{t_{k-1},t_{k}}^{x})-f(\tilde{Y}_{t_{k-1},t_{k}}^{x})]\\
    =&\mathbb{E} \left[f\left(x+\int_{t_{k-1}}^{t_{k}}b(X_{t_{k-1},r}^{x})\,\dif r +A(Z_{t_{k}}-Z_{t_{k-1}})\right)
        - f\left(x+\gamma_{k}b(x) +\frac{\gamma_{k}^{1/\alpha}}{\beta}A\tilde{Z}\right)\right]\\
    =&\mathcal{J}_{1}+\mathcal{J}_{2},
\end{align*}
where
\begin{align*}
    \mathcal{J}_{1}
    :=\mathbb{E} \left[f\left(x+\int_{t_{k-1}}^{t_{k}} b(X_{t_{k-1},r}^{x})\,\dif r+AZ_{\gamma_{k}}\right)
        -f\left(x+\gamma_{k}b(x)+AZ_{\gamma_{k}}\right)\right],
\end{align*}
\begin{align*}
    \mathcal{J}_{2}
    :=& \mathbb{E} \left[f\left(x+\gamma_{k}b(x)+AZ_{\gamma_{k}}\right)-f\left(x+\gamma_{k} b(x)\right)\right]\\
        &-\mathbb{E} \left[f\left(x+\gamma_{k} b(x)+\frac{\gamma_{k}^{1/\alpha}}{\beta} A\tilde{Z}\right)-f\left(x+\gamma_{k} b(x)\right)\right].
\end{align*}
We can bound $\mathcal{J}_{1}$ by using \eqref{2} with $\kappa=1$:
\begin{align*}
    |\mathcal{J}_{1}|
    &\leq\|\nabla f\|_{\infty}\mathbb{E} \left|\int_{t_{k-1}}^{t_{k}}b(X_{t_{k-1},r}^{x})\,\dif r-\gamma_{k} b(x)\right|\\
    &\leq\|\nabla f\|_{\infty}\int_{t_{k-1}}^{t_{k}}\mathbb{E} |b(X_{t_{k-1},r}^{x})-b(x)|\,\dif r\\
    &\leq L\|\nabla f\|_{\infty}\int_{t_{k-1}}^{t_{k}}\mathbb{E} |X_{t_{k-1},r}^{x}-x|\,\dif r\\
    &\leq C(1+|x|)\|\nabla f\|_{\infty}d^{\frac{5}{2\kappa} \vee \frac{3}{2}} \int_{t_{k-1}}^{t_{k}}(r-t_{k-1})^{1/\alpha}\,\dif r \\
    &\leq C(1+|x|)\| \, \nabla f\|_{\infty}d^{\frac{5}{2\kappa} \vee \frac{3}{2}} \gamma_{k}^{1+1/\alpha}.
\end{align*}
For the first term of $\mathcal{J}_{2}$, we use Dynkin's formula (see e.g.\ \cite{CX19}) to get
\begin{align*}
    \mathbb{E} \big[f\big(x+\gamma_{k}b(x)+AZ_{\gamma_{k}}\big)-f\big(x+\gamma_{k} b(x)\big)\big]
    = \int_{0}^{\gamma_{k}}\mathbb{E} \big[\Delta^{\alpha/2,A}f\big(x+\gamma_{k} b(x)+Z_{r}\big)\big]\,\dif r.
\end{align*}
For the second part of $\mathcal{J}_{2}$, observing that $d_{\alpha}= \alpha \sigma_{d-1}^{-1} \beta^{-\alpha}$, Taylor's formula implies that
\begin{align*}
    &\mathbb{E}\left[f\left(x+\gamma_{k} b(x)+\frac{\gamma_{k}^{1/\alpha}}{\beta}A\tilde{Z}\right)-f\big(x+\gamma_{k} b(x)\big)\right]\\
    &\quad = \frac{\gamma_{k}^{1/\alpha}}{\beta}\mathbb{E}\left[\int_{0}^{1}\left\langle\nabla f\left(x+\gamma_{k} b(x)+
    \frac{\gamma_{k}^{1/\alpha}}{\beta}tA\tilde{Z}\right),\: A\tilde{Z}\right\rangle \dif t\right]\\
    &\quad = \frac{\gamma_{k}^{1/\alpha}}{\beta} \int_{|z|\geq1} \int_{0}^{1} \alpha \left\langle\nabla f\big(x+\gamma_{k} b(x)
    +\frac{\gamma_{k}^{1/\alpha}} {\beta}tAz\big),\: Az\right\rangle \frac{\dif t\,\dif z}{\sigma_{d-1} |z|^{\alpha+d}}\\
    &\quad = \frac{\alpha\gamma_{k}}{\sigma_{d-1} \beta^{\alpha}}\int_{|z|\geq\beta^{-1}\gamma_{k}^{1/\alpha}}\int_{0}^{1}
    \left\langle\nabla f\big(x+\gamma_{k} b(x)+tAz\big),\:Az\right\rangle  \frac{\dif t\,\dif z}{|z|^{\alpha+d}}\\
    &\quad = \gamma_{k}\Delta^{\alpha/2,A}f(x+\gamma_{k} b(x))-\mathcal{R},
\end{align*}
where
\begin{align*}
    \mathcal{R}
    := \gamma_{k}d_{\alpha}\int_{|z|<\beta^{-1}\gamma_{k}^{1/\alpha}}\int_{0}^{1}\left\langle\nabla f\big(x+\gamma_{k}b(x)+tAz\big),\: Az\right\rangle
    \frac{\dif t\,\dif z}{|z|^{\alpha+d}}.
\end{align*}
Together, the above estimates yield
\begin{align*}
    |\mathcal{J}_{2}|
    \leq |\mathcal{R}| + \left|\int_{0}^{\gamma_{k}}\mathbb{E} \big[\Delta^{\alpha/2,A} f\big(x+\gamma_{k} b(x)+Z_{r}\big)\big]\,
    \dif r-\gamma_{k}\Delta^{\alpha/2,A}f(x+\gamma_{k} b(x))\right|.
\end{align*}
Further, we have
\begin{align*}
    |\mathcal{R}|
    &= \gamma_{k}d_{\alpha} \left|\int_{|z|<\beta^{-1}\gamma_{k}^{1/\alpha}}\int_{0}^{1} \left\langle\nabla f\big(x+\gamma_{k} b(x)
    +tAz\big)-\nabla f\big(x+\gamma_{k} b(x)\big),\: Az\right\rangle \frac{\dif t\,\dif z}{|z|^{\alpha+d}} \right|\\
    &\leq d_{\alpha} \|A\|_{{\rm op}}\gamma_{k}\int_{|z|<\beta^{-1}\gamma_{k}^{1/\alpha}}\int_{0}^{1} \left|\nabla f\big(x+\gamma_{k} b(x)
    +tAz\big)-\nabla f\big(x+\gamma_{k} b(x)\big)\right| \frac{\dif t\,\dif z}{|z|^{\alpha+d-1}}\\
    &\leq C\|\nabla^{2}f\|_{\mathrm{HS},\infty} d_{\alpha} \|A\|_{{\rm op}}^{2}\gamma_{k}\int_{|z|<\sigma^{-1}\eta^{1/\alpha}} \frac{\dif z}{|z|^{\alpha+d-2}}\\
    &\leq C\|\nabla^{2}f\|_{\mathrm{HS},\infty} d_{\alpha}\sigma_{d-1} \|A\|_{{\rm op}}^{2}\gamma_{k}^{2/\alpha}
    \leq C\|\nabla^{2}f\|_{\mathrm{HS},\infty}\,d^{2}\gamma_{k}^{2/\alpha}.
\end{align*}
By Lemma \ref{laplace}, we also have
\begin{align*}
    &\left|\int_{0}^{\gamma_{k}}\mathbb{E}\left[\Delta^{\alpha/2,A}f\big(x+\gamma_{k} b(x)+Z_{r}\big)\right] \dif r
    - \gamma_{k}\Delta^{\alpha/2,A}f(x+\gamma_{k} b(x))\right|\\
    &\quad \leq\int_{0}^{\gamma_{k}}\mathbb{E}\left|\Delta^{\alpha/2,A}f\big(x+\gamma_{k} b(x)+Z_{r}\big)\big] -\Delta^{\alpha/2,A}f(x+\gamma_{k} b(x))\right| \dif r\\
    &\quad \leq C\|\nabla^{2}f\|_{\mathrm{HS},\infty}d^{2}\int_{0}^{\gamma_{k}}\mathbb{E}\left[|Z_{r}|^{2-\alpha}\right] \dif r\\
    &\quad =C\|\nabla^{2}f\|_{\mathrm{HS},\infty}d^{2}\int_{0}^{\gamma_{k}}\mathbb{E}\left[|Z_{1}|^{2-\alpha}\right] r^{2/\alpha-1}\,\dif r
    \leq C \|\nabla^{2}f\|_{\mathrm{HS},\infty}d^{3-\frac{\alpha}{2}} \,\gamma_{k}^{2/\alpha}.
\end{align*}
Combining all the estimates finishes the proof.
\end{proof}

\subsection{Proof of Lemma  \ref{pro:erg-X}}\label{sec:para}

Recall that $X_0=x$ and
\begin{eqnarray}\label{e:SDE-2}
\dif X_t &=& b(X_t) \dif t + A \dif Z_t,
\end{eqnarray}
and the corresponding generator $\mathcal{L}$ for process $(X_t)_{t\geq 0}$ is defined as
\begin{eqnarray*}
\mathcal{L} f(x) = \langle b(x), \nabla f(x) \rangle + \int_{ \R^d \setminus \{ 0 \}  } [f(x+Az)-f(x)-\langle Az, \nabla f(x) \rangle
1_{\{|z|\leq 1\}}] \frac{d_{\alpha}}{|z|^{d+\alpha}} \dif z,
\end{eqnarray*}
where $\alpha\in(1,2)$. In \cite{Liang2020Gradient}, the drift coefficient $b$ is assumed to satisfy the following condition:
there exist positive constants $K_1,K_2$ and $L_0$, all independent of dimension $d$, such that
\begin{eqnarray}\label{e:b-3}
\langle x-y, b(x)-b(y) \rangle
&\leq&
\left\{
\begin{array}{lll}
K_1 |x-y|^2,  &  \text{if } |x-y|\leq L_0,  \\
-K_2 |x-y|^2, &  \text{if } |x-y| > L_0.
\end{array}
\right.
\end{eqnarray}

In addition, there exists a positive constant $\Lambda \asymp d^{\frac{1}{2}}$ and $\Lambda \geq 1$ such that
\begin{eqnarray*}
\Lambda^{-1}  \ \  \leq \ \  \| A \|_{\rm op} \vee \| A^{-1} \|_{\rm op}  \ \ \leq \ \ \Lambda.
\end{eqnarray*}

We recall some notations in \cite{Liang2020Gradient}. For some small enough constant $\tau \leq \frac{1}{32 \Lambda^5}
 \wedge L_0$ with $L_0$ and $\Lambda$ above, for any $x,y \in \R^d$, define
\begin{eqnarray*}
(x-y)_{\tau} &=& \left(1 \wedge \frac{\tau}{|x-y|} \right) (x-y),
\end{eqnarray*}
and for any $z\in \R^d$, let $\Psi: \R^d \to \R^d$ be a continuous and bijective mapping defined by
\begin{eqnarray*}
\Psi(z) &=& z + A^{-1} (x-y)_{\tau},
\end{eqnarray*}
and then the inverse function $\Psi^{-1}$ is given by
\begin{eqnarray*}
\Psi^{-1}(z) &=& z - A^{-1} (x-y)_{\tau}.
\end{eqnarray*}

Let the measure $\pi_0(\dif z)$ be related to the L\'{e}vy measure $\pi(\dif z)$ and satisfy that $0\leq \pi_0(\dif z) \leq \pi(\dif z)$
and $\int_{\R^d} |z|^2 \pi_0(\dif z)<\infty.$

For any integers $n\geq 1$, define
\begin{eqnarray*}
\mu_{\pi_0,\Psi} \ \ = \ \ \limsup_{n\to \infty} \mu_{\pi_0,n,\Psi}
= \ \ \limsup_{n\to \infty} [\pi_{0,n} \wedge (\pi_{0,n} \Psi)],
\end{eqnarray*}
where $\pi_{0,n}(B)=\int_{B\cap \{|z|>1/n\} }\pi_0(\dif z)$ and $(\pi_{0,n}\Psi)(B)=\pi_{0,n}(\Psi(B))$ for all Borel measurable
set $B\in \mathcal{B}(\R^d)$.

For the solution  $(X_t)_{t\geq 0}$ of SDE \eqref{e:SDE-2} and its operator $\mathcal{L}$, it follows from  \cite[Theorem 2.3 and
Proposition 2.6]{Liang2020Gradient} that there exist the coupling process for $(X_t)_{t\geq 0}$ and the coupling
operator $\widetilde{\mathcal{L}}$ that is defined as follows. For any $h\in \mathcal{C}_b^2(\R^d \times \R^d)$ and $x,y\in \R^d$,
\begin{eqnarray*}
\widetilde{\mathcal{L}} h(x,y)
&=& \langle \nabla_x h(x,y),b(x) \rangle+\langle \nabla_y h(x,y),b(y) \rangle \\
&& +\frac{1}{2} \int_{ \R^d \setminus \{ 0 \} } ( h(x+ Az,y+A \Psi(z))-h(x,y)-\langle \nabla_x h(x,y),A z \rangle 1_{\{ |z|\leq 1  \}  } \\
&& \ \ \ \ \ \ \ \ \ \ \ \ \ \ \ -\langle \nabla_y h(x,y),  A \Psi(z) \rangle 1_{  \{ |\Psi(z)| \leq 1  \}  } ) \mu_{\pi_0,\Psi}(\dif z) \\
&& +\frac{1}{2} \int_{ \R^d \setminus \{ 0 \}  } ( h(x+A z,y+A \Psi^{-1}(z))-h(x,y)-\langle \nabla_x h(x,y),A z \rangle 1_{\{ |z|\leq 1  \}  } \\
&& \ \ \ \ \ \ \ \ \ \ \ \ \ \ \ -\langle \nabla_y h(x,y), A \Psi^{-1}(z) \rangle 1_{  \{ |\Psi^{-1}(z)| \leq 1  \}  } ) \mu_{\pi_0,\Psi^{-1}}(\dif z) \\
&&+ \int_{ \R^d \setminus \{ 0 \} } ( h(x+A z,y+A z)-h(x,y)-\langle \nabla_x h(x,y),A z \rangle 1_{\{ |z|\leq 1  \}  } \\
&& \ \ \ \ \ \ \ \ \ \ \ \ \ \ \ -\langle \nabla_y h(x,y), A z \rangle 1_{  \{ |z| \leq 1  \}  } ) \left[ \pi- \frac{1}{2}\mu_{\pi_0,\Psi} -
\frac{1}{2}\mu_{\pi_0,\Psi^{-1}}  \right](\dif z),
\end{eqnarray*}
where $\nabla_x$ and $\nabla_y$ mean the gradient with respected to variables $x$ and $y$ respectively.

Denote $\pi_0(\dif z) = q(z) \dif z$ and $(\pi_0 \Psi^{-1})(\dif z) = q_{\Psi}(z) \dif z$. And for $r\geq 0$, let
\begin{eqnarray*}
J(r) &=& \inf_{x,y \in \R^d, |x-y|\leq r} \mu_{\pi_0,\Psi}(\R^d)
\ \ = \ \ \inf_{x,y\in \R^d, |x-y|\leq r} \int_{\R^d} [q(z) \wedge q_{\Psi}(z) ] \dif z.
\end{eqnarray*}

For $r\in (0,2L_0]$ and constants $K_1$ and $K_2$ in inequality \eqref{e:b-3}, define
\begin{eqnarray*}
\psi(r) &\leq& \frac{1}{2r} J(\tau \wedge r)(\tau \wedge r)^2, \ \
g_1(r) \ \ = \ \ \int_0^r \frac{1}{\psi(s)}\dif s,
\end{eqnarray*}
and
\begin{eqnarray*}
g(r) &=& \left(1+\frac{2K_1}{(2K_2) \wedge g_1^{-1}(2L_0)} \right)g_1(r).
\end{eqnarray*}

The proof of exponential ergodicity is standard, see more details in \cite{Liang2020A,Liang2020Gradient,Wang2016Lp} and references therein. It follows from \cite[Theorem 4.4]{Liang2020Gradient} that we obtain the following exponential ergodicity for the process
$(X_t)_{t\geq 0}$ under the Wasserstein-1 distance.
\begin{lemma}\label{thm:ee-para}
Let the process $(X_t)_{t\geq 0}$ be in \eqref{e:SDE-2}, under \eqref{e:b-3}, there exist positive constants $c_{1}$ and
$c_{2}$ independent of $d$ such that for any $x,y\in \R^d$ and $t>0$
\begin{eqnarray*}
W_1( \mathcal{L}(X_t^x), \mathcal{L}(X_t^y) )
&\leq&c_{1}\exp( c_{2}d^{4-2\alpha}2^{d}) e^{-c_{1}d^{2\alpha-4} 2^{-d} \exp(-c_{2}d^{4-2\alpha}2^d)  t} |x-y|.
\end{eqnarray*}
\end{lemma}

\begin{proof} Recall that $\pi(\dif z) = \frac{d_{\alpha}}{|z|^{d+\alpha}} \dif z$, we choose
\begin{eqnarray*}
\pi_0 (\dif z) &=& 1_{\{-\frac{1}{8} < z_1 < \frac{1}{8}, |z|\leq 1 \}} \frac{d_{\alpha}}{|z|^{d+\alpha}} \dif z,
\end{eqnarray*}
where $z_1$ is the first element of vector $z$, that is, $z=(z_1,z_2,\cdots,z_d)$. Thus, we know
\begin{eqnarray*}
q(z) &=& 1_{\{-\frac{1}{8} < z_1 < \frac{1}{8}, |z|\leq 1 \}} \frac{d_{\alpha}}{|z|^{d+\alpha}}, \\
q_{\Psi}(z) &=& 1_{\{-\frac{1}{8} < (\Psi^{-1}(z))_1 <  \frac{1}{8}, |\Psi^{-1}(z)|\leq 1 \}} \frac{d_{\alpha}}{|\Psi^{-1}(z)|^{d+\alpha}},
\end{eqnarray*}
satisfy
\begin{eqnarray*}
q(z) \wedge q_{\Psi}(z) &\geq & d_{\alpha} 1_{\{-\frac{1}{8}< z_1< \frac{1}{8}, -\frac{1}{8} < (\Psi^{-1}(z))_1 < \frac{1}{8}, |z|\leq 1,
|\Psi^{-1}(z)|\leq 1 \}} \left[\frac{1}{|z|^{d+\alpha}} \wedge \frac{1}{|\Psi^{-1}(z)|^{d+\alpha}} \right].
\end{eqnarray*}

We give the estimate for function $J(r)$ firstly. For $x,y\in \R^d$ with $|\tl{x}-\tl{y}| \leq \tau$, then

(i) if $(A^{-1}(\tl{x}-\tl{y}))_1 \leq 0$, with some calculations, one gets that
\begin{eqnarray*}
q(z) \wedge q_{\Psi}(z) &\geq & \frac{d_{\alpha}}{2^{d+\alpha}} 1_{\{z_1>0, \Lambda|\tl{x}-\tl{y}|\leq |z| \leq \frac{1}
{16\Lambda^2} \} } \frac{1}{|z|^{d+\alpha}},
\end{eqnarray*}
that is,
\begin{eqnarray*}
\int_{\R^d} [q(z) \wedge q_{\Psi}(z)] \dif z
&\geq& \frac{d_{\alpha} \sigma_{d-1}}{\alpha 2^{d+\alpha+1}} \left[ \frac{1}{(\Lambda|\tl{x}-\tl{y}|)^{\alpha}} - \frac{1}
{(\frac{1}{16\Lambda^2})^{\alpha}}\right] \\
&\geq& \frac{d_{\alpha} \sigma_{d-1}}{\alpha 2^{d+\alpha+1}} (\frac{1}{\Lambda^{\alpha}}-\frac{1}{2^{\alpha}\Lambda^{3\alpha}}) |\tl{x}-\tl{y}|^{-\alpha},
\end{eqnarray*}
where the last inequality holds provided $|\tl{x}-\tl{y}| \leq \tau \leq \frac{1}{32 \Lambda^5}$.

(ii) if $(A^{-1}(\tl{x}-\tl{y}))_1>0$, then
\begin{eqnarray*}
q(z) \wedge q_{\Psi}(z) &\geq & \frac{d_{\alpha}}{2^{d+\alpha}} 1_{\{ (\Psi^{-1}(z))_1>0, \Lambda|\tl{x}-\tl{y}|\leq |\Psi^{-1}(z)|
\leq \frac{1}{16\Lambda^2} \} } \frac{1}{|\Psi^{-1}(z)|^{d+\alpha}},
\end{eqnarray*}
that is,
\begin{eqnarray*}
\int_{\R^d} [q(z) \wedge q_{\Psi}(z)] \dif z
&\geq& \frac{d_{\alpha} \sigma_{d-1}}{\alpha 2^{d+\alpha+1}} \left[ \frac{1}{(\Lambda|\tl{x}-\tl{y}|)^{\alpha}} -
 \frac{1}{(\frac{1}{16\Lambda^2})^{\alpha}}\right] \\
&\geq& \frac{d_{\alpha} \sigma_{d-1}}{\alpha 2^{d+\alpha+1}} (\frac{1}{\Lambda^{\alpha}}-\frac{1}{2^{\alpha}
\Lambda^{3\alpha}}) |\tl{x}-\tl{y}|^{-\alpha},
\end{eqnarray*}
where the last inequality holds provided $|\tl{x}-\tl{y}| \leq \tau \leq \frac{1}{32\Lambda^5}$.

Combining above two cases, for any $0<r\leq \tau$, we know
\begin{eqnarray*}
\int_{\R^d} [q(z) \wedge q_{\Psi}(z)] \dif z
&\geq& \frac{d_{\alpha} \sigma_{d-1}}{\alpha 2^{d+\alpha+1}} (\frac{1}{\Lambda^{\alpha}}-\frac{1}{2^{\alpha}\Lambda^{3\alpha}}) |\tl{x}-\tl{y}|^{-\alpha}
\ \ = \ \ \widehat{C}_1 |\tl{x}-\tl{y}|^{-\alpha}
\end{eqnarray*}
with $\widehat{C}_1=\frac{d_{\alpha} \sigma_{d-1}}{\alpha 2^{d+\alpha+1}} (\frac{1}{\Lambda^{\alpha}}-\frac{1}{2^{\alpha}
\Lambda^{3\alpha}})$. Thus, for any $0<r\leq \tau$, we know
\begin{eqnarray*}
J(r) &=&\inf_{\tl{x},\tl{y}\in \R^d,|\tl{x}-\tl{y}|\leq r} \int_{\R^d} [q(z) \wedge q_{\Psi}(z)] \dif z
\ \ \geq  \ \ \widehat{C}_1 r^{-\alpha}.
\end{eqnarray*}
It implies that we can choose function $\psi$ as follows:
\begin{eqnarray*}
\psi(r)
&=&
\left\{
\begin{array}{lll}
\frac{\widehat{C}_1}{2} r^{1-\alpha},  &  \text{if } r\in (0,\tau],  \\
\frac{\widehat{C}_1}{2}\tau^{2-\alpha} r^{-1}, &  \text{if } r\in (\tau,2L_0].
\end{array}
\right.
\end{eqnarray*}

Then, we know
\begin{eqnarray*}
g_1(2L_0) &=& \left(\int_0^{\tau} + \int_{\tau}^{2L_0}  \right) \frac{1}{\psi(s)} \dif s
\ \ = \ \ 2\alpha^{-1} \widehat{C}_1^{-1} \tau^{\alpha}+ \tau^{\alpha-2} \widehat{C}_1^{-1} (4L_0^2 - \tau^2),  \\
g(2L_0) &=& \left(1+\frac{2K_1}{(2K_2) \wedge g_1^{-1}(2L_0)} \right)g_1(2L_0).
\end{eqnarray*}

Using \cite[Theorem 4.4]{Liang2020Gradient}, we obtain that
\begin{eqnarray*}
W_1( \mathcal{L}(X_t^x), \mathcal{L}(X_t^y) ) &\leq& \widehat{C} e^{-\lambda t} |x-y|,
\end{eqnarray*}
where
\begin{eqnarray*}
\widehat{C} &=& \frac{1}{2}\left( 1+\exp\left\{ g(2L_0) [ (2K_2) \wedge g_1^{-1}(2L_0) ] \right\} \right) ,
\end{eqnarray*}
and
\begin{eqnarray*}
\lambda &=& \frac{(2K_2) \wedge g_1^{-1}(2L_0)}{1+\exp\left\{ g(2L_0)[ (2K_2) \wedge g_1^{-1}(2L_0) ] \right\}   }.
\end{eqnarray*}

Recall that  $\widehat{C}_1=\frac{d_{\alpha} \sigma_{d-1}}{\alpha 2^{d+\alpha+1}} (\frac{1}{\Lambda^{\alpha}}-\frac{1}{2^{\alpha}
\Lambda^{3\alpha}})$, $d_{\alpha}= 2^{\alpha}\Gamma(\frac{d+\alpha}{2})\pi^{-\frac{d}{2}}|\Gamma(-\frac{\alpha}{2})|^{-1}$
and $\sigma_{d-1}=2\pi^{\frac{d}{2}} \Gamma^{-1}(\frac{d}{2})$, then there exists a positive constant $C$ independent of $d$ such that
\begin{eqnarray*}
\widehat{C}_1 &\leq& Cd^{1-\frac{\alpha}{2}} 2^{-d}.
\end{eqnarray*}
Moreover, there exist positive constants $c_{1}$ and $c_{2}$ both independent of $d$  such that
\begin{eqnarray*}
\widehat{C} \ \ \leq \ \  c_{1}\exp\left( c_{2}d^{4-2\alpha }2^{d} \right), \qquad
\lambda \ \ \geq \ \ c_{1}d^{2\alpha-4}2^{-d} \exp\left(-c_{2}d^{4-2\alpha} 2^d \right).
\end{eqnarray*}
The proof is complete.
\end{proof}

Now, we are in the position to prove Lemma \ref{pro:erg-X}.

\begin{proof}[Proof of Lemma \ref{pro:erg-X}]
In \eqref{e:b-3}, let $K_{1}= L$, $K_{2}=\frac{\theta_{1}}{2}$ and $L_{0}=\sqrt{\frac{2K}{\theta_{1}}}$, then the
desired result following from {\bf Assumption} \ref{assump-1} and Lemma \ref{thm:ee-para}.
\end{proof}

\subsection{Proof of Lemma \ref{lem:cal-ss}}

For $n\geq1$, set $u_{n}=\frac{v_{n}}{\gamma_{n}^{\theta}}$. Then
\begin{align}\label{Eq:un}
u_{n+1}=u_{n}\lambda_n+\gamma_{n+1},
\end{align}
where $\lambda_n=\left(\frac{\gamma_{n}}{\gamma_{n+1}}\right)^{\theta}e^{-\rho\gamma_{n+1}}$. Recall that
$\omega = \limsup\limits_{k\to \infty} \frac{\gamma^{\theta}_k-\gamma^{\theta}_{k+1}}{\gamma_{k+1}^{1+\theta}}
< +\infty$ and $\rho>\omega$.  Then for any constant $c\in(\omega,\rho)$, there exists an $n_{0}\in\mathbb{N}$
such that for all $n\geq n_{0}$,
\begin{eqnarray}\label{e:claim-ss1}
\frac{\gamma^{\theta}_n}{\gamma^{\theta}_{n+1}} = 1+\frac{\gamma^{\theta}_{n} - \gamma^{\theta}_{n+1}}
{\gamma^{1+\theta}_{n+1}}  \gamma_{n+1}
\  \leq  \ 1+ c \gamma_{n+1}
\  \leq \ e^{c \gamma_{n+1}}.
\end{eqnarray}
Hence,
\begin{eqnarray*}
\lambda_n \leq e^{-(\rho-c)\gamma_{n+1}}.
\end{eqnarray*}
This and (\ref{Eq:un}) imply
\begin{eqnarray*}
u_{n+1} \leq e^{-(\rho-c)\gamma_{n+1}}u_{n}+\gamma_{n+1}.
\end{eqnarray*}
Consequently,
\begin{eqnarray*}
e^{(\rho-c)t_{n+1}} u_{n+1}
\leq  e^{(\rho-c)t_{n}}u_n + M e^{(\rho-c)t_{n}}\gamma_{n+1},
\end{eqnarray*}
where the constant $M=\sup_{k\geq1}e^{(\rho-c)\gamma_{k}}$. Hence, by induction, we obtain
\begin{align*}
e^{(\rho-c)t_{n}} u_{n}\leq e^{(\rho-c)t_{n_{0}}} u_{n_{0}}+M\int_{t_{n_{0}}}^{t_{n}}e^{(\rho-c)u}du\leq e^{(\rho-c)t_{n_{0}}}
u_{n_{0}}+\frac{M}{\rho-c}e^{(\rho-c)t_{n}}.
\end{align*}
Consequently,
\begin{align} \label{Eq:un2}
 u_{n}\leq e^{(\rho-c)(t_{n_{0}} - t_n)} u_{n_{0}}   +\frac{M}{\rho-c}.
\end{align}
We take $c=\frac{\rho+\omega}{2}$, then $M=e^{(\rho-c)\gamma_1} = e^{\frac{\rho-\omega} 2 \gamma_1}$.
Letting  $n \to \infty$ in (\ref{Eq:un2}) yields   (\ref{convergence1}).

For proving (\ref{convergence2}), we set $\omega_{n}=\frac{e^{-\rho t_{n}}}{\gamma_{n}^{\theta}}$. Then, by (\ref{e:claim-ss1}),
we see that for $n\geq n_{0}$,
\begin{align} \label{Eq:un3}
\omega_{n+1}=\omega_{n}e^{-\rho\gamma_{n+1}}
\left(\frac{\gamma_{n}}{\gamma_{n+1}}\right)^{\theta}
\leq \omega_{n}e^{(c-\rho)\gamma_{n+1}}
\leq\omega_{n_{0}}e^{(c-\rho)(t_{n+1}-t_{n_{0}})}.
\end{align}
Hence (\ref{convergence2})  follows from (\ref{Eq:un3}) and the facts that $\sum_{k=n_{0}}^{\infty}\gamma_{k}=\infty$ and $c < \rho$.

Finally, we prove the last inequality \eqref{convergence3}. By the same argument as the proof of \eqref{e:claim-ss1}, we have
$$\sup_{n\geq1}\frac{\gamma_{n}^{\frac{1}{\alpha}}}{\gamma_{n+1}^{\frac{1}{\alpha}}}<1+c\gamma_{1}.$$
This implies that for any $i\leq n-1$,
\begin{align*}
\frac{(t_{n}-t_{i-1})^{\frac{1}{\alpha}}}{(t_{n}-t_{i})^{\frac{1}{\alpha}}}=\frac{(t_{n}-t_{i}+\gamma_{i})^{\frac{1}{\alpha}}}{(t_{n}-t_{i})^{\frac{1}{\alpha}}}
\leq 2+\frac{2\gamma_{i}^{\frac{1}{\alpha}}}{(t_{n}-t_{i})^{\frac{1}{\alpha}}}\leq2+\frac{2\gamma_{i}^{\frac{1}{\alpha}}}
{\gamma_{i+1}^{\frac{1}{\alpha}}}\leq C.
\end{align*}
Hence, we have
\begin{align*}
\sum_{i=n^*+1}^{n-1}(t_n-t_i)^{-\frac{1}{\alpha}} \gamma_i^{1+\theta}
\leq&C\sum_{i=n^*+1}^{n-1}(t_n-t_{i-1})^{-\frac{1}{\alpha}} \gamma_i^{1+\theta}\\
\leq&C\gamma_{n^{*}}^{\theta}\int^{t_{n-1}}_{t_{n^{*}}}(t_{n}-t)^{-\frac{1}{\alpha}}\dif t\leq C\gamma_{n^{*}}^{\theta}.
\end{align*}
In addition, \eqref{e:claim-ss1} implies
\begin{align*}
\frac{\gamma_{n^{*}}^{\theta}}{\gamma_{n}^{\theta}}
=\prod_{i=n^{*}}^{n-1}(\frac{\gamma_{i}}{\gamma_{i+1}})^{\theta}
\leq e^{c(t_{n}-t_{n^{*}})}
\leq e^{\rho(1+\gamma_{1})}.
\end{align*}
Therefore, we have
\begin{align*}
\sum_{i=n^*+1}^{n-1}(t_n-t_i)^{-\frac{1}{\alpha}} \gamma_i^{1+\theta}
\leq Ce^{\rho(1+\gamma_{1})}
\gamma_{n}^{\theta}.
\end{align*}
This proves \eqref{convergence3}.
\qed
~\\

\begin{appendix}
\section{Lower bound for the stable OU process}\label{seclow}

We consider the stable Ornstein-Uhlenbeck (OU) process defined by
\begin{eqnarray}\label{e:OU}
\dif X_t = -X_t \dif t + \dif Z_t, \qquad X_0=x,
\end{eqnarray}
where $(Z_t)_{t\geq 0}$ is a one-dimensional symmetric $\alpha$-stable process with $\alpha\in(1,2)$.
It is known that the solution to \eqref{e:OU} is
\begin{eqnarray}\label{e:OU-solu}
X_t^x = e^{-t}x + e^{-t} \int_0^t e^s \dif Z_s, \quad \forall t\geq 0.
\end{eqnarray}

The EM scheme with the Pareto distribution for $(X_t)_{t\geq 0}$ in \eqref{e:OU} is given as follows:
$\tilde{Y}_0=x$ and for any integers $n\geq 1$, one has
\begin{align}\label{e-OU1}
\tilde{Y}_{t_{n+1}}=\tilde{Y}_{t_{n}}-\gamma_{n+1}\tilde{Y}_{t_{n}}
+\frac{\gamma_{n+1}^{\frac{1}{\alpha}}}{\beta}\tl{Z}_{n+1},
\end{align}
where $\beta=\left(\frac{\alpha}{2d_{\alpha}}\right)^{1/\alpha}$ and $(\tilde{Z}_{n})_{n\geq 1}$ is a
sequence of i.i.d. random variables with density
\begin{align*}
p(z)=\frac{\alpha}{2|z|^{\alpha+1}}{\bf 1}_{(1,\infty)}(|z|).
\end{align*}
Then, we have the following lemma.

\begin{lemma}\label{lower bound}
Denote the invariant measure of $(X_{t})_{t\geq0}$ defined in \eqref{e:OU} by $\nu$ and let $(\tl{Y}_{t_{n}})_{n\in\mathbb{N}_{0}}$
be defined by \eqref{e-OU1}. Suppose Assumption \ref{stepsize} holds with $\theta=\frac{2-\alpha}{\alpha}$ and $\omega<\rho\wedge1$.
Then, for any bounded $x\in\mathbb{R}$ and $\alpha\in(1,2)$, we have
\begin{align*}
0<\liminf_{n\rightarrow\infty}\frac{W_{1}(\nu,\tilde{Y}_{t_{n}}^{x})}{\gamma_{n}^{\frac{2}{\alpha}-1}}
\leq\limsup_{n\rightarrow\infty}\frac{W_{1}(\nu,\tilde{Y}_{t_{n}}^{x})}{\gamma_{n}^{\frac{2}{\alpha}-1}}
<\infty.
\end{align*}
\end{lemma}

\begin{proof}
From \eqref{e:OU-solu}, it is easy to  see that the invariant measure $\nu$ has the same distribution as the random variable
$\alpha^{-\frac{1}{\alpha}}Z_{1}$, and from \eqref{e-OU1}, we can derive that
\begin{align*}
\tilde{Y}_{t_{n}}^{x}=&\prod_{j=1}^{n}(1-\gamma_{j})x+\sum_{j=0}^{n-1}\frac{\gamma_{n-j}^{\frac{1}{\alpha}}\tilde{Z}_{n-j}}
{\beta}\prod_{k=n+1-j}^{n}(1-\gamma_{k})\\
=&\prod_{j=1}^{n}(1-\gamma_{j})x+\sum_{j=1}^{n}\frac{\gamma_{j}^{\frac{1}{\alpha}}\tilde{Z}_{j}}{\beta}\prod_{k=j+1}^{n}(1-\gamma_{k}).
\end{align*}
Denote $\phi(\lambda)=\mathbb{E}[e^{i \lambda \tilde{Z}_{1}}]$, then we have
\begin{align*}
\mathbb{E}[e^{i\lambda\tilde{Y}_{t_{n}}^{x}}]
=&e^{i\lambda\prod_{j=1}^{n}(1-\gamma_{j})x}\prod_{j=1}^{n}\mathbb{E}\left[e^{i\lambda\frac{\gamma_{j}^{\frac{1}{\alpha}}\tilde{Z}_{j}}
{\beta}\prod_{k=j+1}^{n}(1-\gamma_{k})}\right]\nonumber\\
=&e^{i\lambda\prod_{j=1}^{n}(1-\gamma_{j})x}\prod_{j=1}^{n}\phi\left(\frac{\gamma_{j}^{\frac{1}{\alpha}}}{\beta}\prod_{k=j+1}^{n}
(1-\gamma_{k})\lambda\right).
\end{align*}
By \cite[(B.4)]{Chen2022stableEM}, for $|\lambda|\leq 1$ and $0<\gamma_{1}<1\wedge\beta^{\alpha}$, we have
\begin{align*}
&\log\phi\left(\frac{\gamma_{j}^{\frac{1}{\alpha}}}{\beta}\prod_{k=j+1}^{n}(1-\gamma_{k})\lambda\right)\\
\geq&\log\left(1-\gamma_{j}\prod_{k=j+1}^{n}(1-\gamma_{k})^{\alpha}|\lambda|^{\alpha}+\frac{c\alpha}{(2-\alpha)\beta^{2}}
\gamma_{j}^{\frac{2}{\alpha}}\prod_{k=j+1}^{n}(1-\gamma_{k})^{2}|\lambda|^{2}\right),
\end{align*}
where $c=\inf_{0<u\leq1}\frac{1-\cos u}{u^{2}}>0$. Since
\begin{align*}
\lim_{x\downarrow0}\frac{\log\left(1-x+\frac{c\alpha}{(2-\alpha)\beta^{2}}x^{\frac{2}{\alpha}}\right)+x}
{x^{\frac{2}{\alpha}}}=\frac{c\alpha}{(2-\alpha)\beta^{2}},
\end{align*}
there is a positive constant $C=C(\alpha)$ such that for small enough $x>0$
\begin{align*}
\log\left(1-x+\frac{c\alpha}{(2-\alpha)\beta^{2}}x^{\frac{2}{\alpha}}\right)\geq -x+Cx^{\frac{2}{\alpha}}.
\end{align*}
These imply that for all $|\lambda|\leq 1$ and large enough $n$ or $j$, we have
\begin{align*}
\log\phi\left(\frac{\gamma_{j}^{\frac{1}{\alpha}}}{\beta}\prod_{k=j+1}^{n}(1-\gamma_{k})\lambda\right)
\geq-\gamma_{j}\prod_{k=j+1}^{n}(1-\gamma_{k})^{\alpha}|\lambda|^{\alpha}+C\gamma_{j}^{\frac{2}{\alpha}}
\prod_{k=j+1}^{n}(1-\gamma_{k})^{2}|\lambda|^{2}.
\end{align*}
Hence, we have
\begin{align}\label{com1}
\log\mathbb{E}[e^{i\lambda\left(\tilde{Y}_{t_{n}}^{x}-x\prod_{j=1}^{n}(1-\gamma_{j})\right)}]
=&\sum_{j=1}^{n}\log\phi\left(\frac{\gamma_{j}^{\frac{1}{\alpha}}}{\beta}\prod_{k=j+1}^{n}(1-\gamma_{k})\lambda\right)\nonumber\\
\geq&-\sum_{j=1}^{n}\gamma_{j}\prod_{k=j+1}^{n}(1-\gamma_{k})^{\alpha}|\lambda|^{\alpha}
+C\sum_{j=1}^{n}\gamma_{j}^{\frac{2}{\alpha}}\prod_{k=j+1}^{n}(1-\gamma_{k})^{2}|\lambda|^{2}.
\end{align}

Now, we analyze the last two items in \eqref{com1}. For the first term, by the inequality $1-x \le e^{-x}$,
we have
\begin{align} \label{com3}
\sum_{j=1}^{n}\gamma_{j}\prod_{k=j+1}^{n}
(1-\gamma_{k})^{\alpha}
\le \sum_{j=1}^{n}\gamma_{j} e^{-\alpha(t_n-t_j)}  \le \frac 1{\alpha}(1-e^{-\alpha t_n})
 +O(\gamma_n),
\end{align}
where the last inequality is obtained by
\begin{eqnarray*}
&& \sum_{j=1}^{n}\gamma_{j} e^{-\alpha(t_n-t_j)}- \frac{1}{\alpha}+\frac{1}{\alpha} e^{-\alpha t_n}
\ \ = \ \ e^{-\alpha t_n} \left[ \sum_{j=1}^{n}\gamma_{j} e^{\alpha t_j}- \int_0^{t_n} e^{\alpha s} \dif s \right] \\
&=& e^{-\alpha t_n} \sum_{j=1}^{n} \int_{t_{j-1}}^{t_j} (e^{\alpha t_j}
- e^{\alpha s}) \dif s
\ \ \leq \ \  e^{-\alpha t_n} \sum_{j=1}^{n} \int_{t_{j-1}}^{t_j} \alpha e^{\alpha t_j} (t_j-s) \dif s  \\
&=& \frac{\alpha}{2} e^{-\alpha t_n} \sum_{j=1}^{n}  e^{\alpha t_j} \gamma_j^2
\ \ \leq \ \ C \gamma_n,
\end{eqnarray*}
where the first inequality holds from the mean value theorem and the last inequality holds from \eqref{convergence1}.

Now we estimate the second term in \eqref{com1}.  Notice that, for some $x_0>0$, we have $1-x\geq e^{-\frac{3}{2}x}$ for
all $x \in (0, x_0)$. Let $j_0=\inf\{j>0: \gamma_j<x_0\}$. Then
\begin{equation}
\begin{split}\label{com4}
\sum_{j=1}^{n} \gamma_{j}^{\frac{2}{\alpha}}\prod_{k=j+1}^{n}
(1-\gamma_{k})^{2}\gamma_{n}^{-\frac{2-\alpha}{\alpha}}
\geq & \sum_{j=j_0}^{n}\gamma_{j} \prod_{k=j+1}^{n}
(1-\gamma_{k})^{2} \left(\frac{\gamma_j}{\gamma_{n}}\right)^{\frac{2-\alpha}{\alpha}} \\
 \ge & \sum_{j=j_0}^{n}\gamma_{j} \prod_{k=j+1}^{n} (1-\gamma_{k})^{2}   \\
\geq & \sum_{j=j_0}^{n} \gamma_{j} e^{-3(t_n-t_j)} \\
\ge & \int_{t_{j_0}}^{t_n} e^{-3(t_n-s)} \dif s= \frac 1 3 \big(1 -  e^{-3 (t_n - t_{j_0})}\big).
\end{split}
\end{equation}
%where $k_0=\gamma_1+...+\gamma_{j_0}=O(1)$ since $\gamma_n \rightarrow 0$ and $t_n \rightarrow \infty$ as $n \rightarrow \infty$.
% \footnote{\textcolor{blue}{Is it easy to see (e.g., monotonicity) why the two limits in the first line exist?}}

Combining \eqref{com1}, \eqref{com3}, \eqref{com4} and the easy fact $e^{-\alpha t_n}=o(\gamma^{\frac{2-\alpha}{\alpha}}_n)$,
we have%\footnote{\textcolor{blue}{Is it obvious that the characteristic function
%$\mathbb{E}[e^{i\lambda\tilde{Y}_{t_{n}}}]$ is positive? Should we add an explanation?}}
\begin{align*}
&\log\mathbb{E}[e^{i\lambda\left(\tilde{Y}_{t_{n}}^{x}-x\prod_{j=1}^{n}(1-\gamma_{j})\right)}]\geq-\frac{1}{\alpha}|\lambda|^{\alpha}
+|\lambda|^{\alpha}o\left(\gamma_{n}^{\frac{2-\alpha}{\alpha}}\right)
+|\lambda|^{2}O\left(\gamma_{n}^{\frac{2-\alpha}{\alpha}}\right),
\end{align*}
where $f(n)=o(g(n))$ as $n\rightarrow\infty$ means that $\lim_{n\rightarrow\infty}\frac{f(n)}{g(n)}=0$, and $f(n)=O(g(n))$ as $n\rightarrow\infty$ means that $\lim_{n\rightarrow\infty}\frac{f(n)}{g(n)}<\infty$. By the inequality $e^{x}\geq1+x$ for any
 $x\in\mathbb{R}$, we further have
\begin{align*}
\mathbb{E}[e^{i\lambda\left(\tilde{Y}_{t_{n}}^{x}-x\prod_{j=1}^{n}(1-\gamma_{j})\right)}]\geq&\exp\left[-\frac{1}{\alpha}|\lambda|^{\alpha}
+|\lambda|^{\alpha}o\left(\gamma_{n}^{\frac{2-\alpha}{\alpha}}\right)
+|\lambda|^{2}O\left(\gamma_{n}^{\frac{2-\alpha}{\alpha}}\right)\right]\\
\geq&e^{-\frac{|\lambda|^{\alpha}}{\alpha}}\left[1+o\left(\gamma_{n}^{\frac{2-\alpha}{\alpha}}\right)+|\lambda|^{2-\alpha}
O\left(\gamma_{n}^{\frac{2-\alpha}{\alpha}}\right)\right].
\end{align*}
In addition, since
\begin{align*}
\prod_{j=1}^{n}(1-\gamma_{j})=1-\sum_{j=1}^{n}\gamma_{j}\prod_{k=j+1}^{n}(1-\gamma_{k})
\end{align*}
and
\begin{align*}
\lim_{n\rightarrow\infty}\left(1-\sum_{j=1}^{n}\gamma_{j}\prod_{k=j+1}^{n}(1-\gamma_{k})\right)e^{t_{n}}
=&\lim_{n\rightarrow\infty}\left(1-\int_{0}^{t_{n}}e^{-(t_{n}-s)}\dif s\right)e^{t_{n}}\\
=&\lim_{n\rightarrow\infty}\left(1-\int_{0}^{t_{n}}e^{-s}\dif s\right)e^{t_{n}}=1,
\end{align*}
one can derive from \eqref{convergence2} that
\begin{align*}
\lim_{n\rightarrow\infty}\frac{\prod_{j=1}^{n}(1-\gamma_{j})}{\gamma_{n}^{\frac{2-\alpha}{\gamma_{n}}}}
=\lim_{n\rightarrow\infty}\frac{e^{-t_{n}}}{\gamma_{n}^{\frac{2-\alpha}{\alpha}}}=0,
\end{align*}
which implies that for large enough $n$,
\begin{align}\label{com2}
\prod_{j=1}^{n}(1-\gamma_{j})=o\left(\gamma_{n}^{\frac{2-\alpha}{\alpha}}\right).
\end{align}
Hence, for large enough $n$ and any bounded $x$, we have
\begin{align*}
&\int_{-1}^{1}\left(\mathbb{E}[e^{i\lambda
\left(\alpha^{-\frac{1}{\alpha}}Z_{1}
-x\prod_{j=1}^{n}(1-\gamma_{j})\right)}]
-\mathbb{E}[e^{i\lambda\alpha^{-\frac{1}{\alpha}}
Z_{1}}]\right)\dif \lambda\\
=&\int_{-1}^{1}\mathbb{E}[e^{i\lambda\alpha^{
-\frac{1}{\alpha}}Z_{1}}]\left[\cos\left(\lambda x\prod_{j=1}^{n}(1-\gamma_{j})\right)-1\right]\dif \lambda\\
=&
|x|^2 O\left(\prod_{j=1}^{n}(1-\gamma_{j})^2\right)
\int_{-1}^{1}\mathbb{E}[e^{i\lambda\alpha^{
-\frac{1}{\alpha}}Z_{1}}]|\lambda|^2 \dif\lambda
=o\left(\gamma_{n}^{\frac{2-\alpha}{\alpha}}\right),
\end{align*}
where the second equality is by the fact
\begin{align*}
\lim_{r\downarrow0}\frac{\cos r-1}{r^2}=-\frac{1}{2}.
\end{align*}
Thus we have
\begin{align}\label{Ageq}
&\int_{-1}^{1}\left(\mathbb{E}[e^{i\lambda\left(\tilde{Y}_{t_{n}}^{x}-x\prod_{j=1}^{n}(1-\gamma_{j})\right)}]
-\mathbb{E}[e^{i\lambda\left(\alpha^{-\frac{1}{\alpha}}Z_{1}-x\prod_{j=1}^{n}(1-\gamma_{j})\right)}]\right)\dif \lambda\nonumber\\
=&\int_{-1}^{1}\left(\mathbb{E}[e^{i\lambda\left(\tilde{Y}_{t_{n}}^{x}-x\prod_{j=1}^{n}(1-\gamma_{j})\right)}]
-\mathbb{E}[e^{i\lambda\alpha^{-\frac{1}{\alpha}}Z_{1}}]\right)\dif \lambda\nonumber\\
&+\int_{-1}^{1}\left(\mathbb{E}[e^{i\lambda\left(\alpha^{-\frac{1}{\alpha}}Z_{1}-x\prod_{j=1}^{n}(1-\gamma_{j})\right)}]
-\mathbb{E}[e^{i\lambda\alpha^{-\frac{1}{\alpha}}Z_{1}}]\right)\dif \lambda\nonumber\\
\geq&\int_{-1}^{1}e^{-\frac{1}{\alpha}|\lambda|^{\alpha}}\left[o\left(\gamma_{n}^{\frac{2-\alpha}{\alpha}}\right)
+|\lambda|^{2-\alpha}O\left(\gamma_{n}^{\frac{2-\alpha}{\alpha}}\right)\right]\dif\lambda+o\left(\gamma_{n}^{\frac{2-\alpha}{\alpha}}\right)\nonumber\\
=&o\left(\gamma_{n}^{\frac{2-\alpha}{\alpha}}\right)+O\left(\gamma_{n}^{\frac{2-\alpha}{\alpha}}\right)
=O\left(\gamma_{n}^{\frac{2-\alpha}{\alpha}}\right).
\end{align}
Then, by the same argument as in the proof of \cite[Proposition B.1]{Chen2022stableEM} and taking
\begin{align*}
h(y)=\frac{1}{M}\left(\frac{\sin\left(y-x\prod_{j=1}^{n}(1-\gamma_{j})\right)}{y-x\prod_{j=1}^{n}(1-\gamma_{j})}
{\bf 1}_{\{y\neq x\prod_{j=1}^{n}(1-\gamma_{j})\}}+{\bf 1}_{\{y=x\prod_{j=1}^{n}(1-\gamma_{j})\}}\right),\quad y\in\mathbb{R},
\end{align*}
with $M=\sup_{z\in\mathbb{R}\backslash\{0\}}|\frac{z\cos z-\sin z}{z^{2}}|\in(0,\infty)$, the desired result follows.
\end{proof}
\end{appendix}

\section*{Acknowledgements}
P. Chen is supported by NSFC grant (NO. 12301176) and NSF of Jiangsu Province of China grant (NO. BK20220867).
X. Jin was supported in part by the Fundamental Research Funds for the Central Universities grants (JZ2022HGQA0148;
JZ2023HGTA0170). Y. Xiao is supported in part by the NSF grant DMS-2153846. L. Xu is supported by NNSFC grant
(No. 12071499) and University of Macau grant MYRG2018-00133-FST.

\end{document}